\documentclass[12pt]{article}
\usepackage{amsmath}
\usepackage{amssymb}
\usepackage{amscd}
\begin{document}

\newtheorem{theorem}{Theorem}[section]
\newtheorem{remark}[theorem]{Remark}
\newtheorem{mtheorem}[theorem]{Main Theorem}
\newtheorem{observation}[theorem]{Observation}
\newtheorem{proposition}[theorem]{Proposition}
\newtheorem{lemma}[theorem]{Lemma}
\newtheorem{reduction-lemma}[theorem]{Reduction-Lemma}
\newtheorem{step-lemma}[theorem]{Step-Lemma}

\newtheorem{note}[theorem]{}
\newtheorem{extlemma}[theorem]{Ext-Lemma}
\newtheorem{corollary}[theorem]{Corollary}
\newtheorem{example}[theorem]{Example}
\newtheorem{definition}[theorem]{Definition}

\renewcommand{\labelenumi}{(\roman{enumi})}
\newcommand{\dach}[1]{\hat{\vphantom{#1}}}
\numberwithin{equation}{section}

\def\Z{{ \mathbb Z}}
\def\N{{ \mathbb N}}
\def\BB{ \mathbb B}
\def\DD{ \mathbb D}
\def\GG{ \mathbb G}
\def\BBB{ B_\BB}
\def\R{{\bf R}}
\def\D{{\hat{D}}}
\def\Q{{\mathbb Q}}
\def\G{\hat{G}}
\def\C{\hat{C}}
\def\T{{\cal T}}
\def\V{{\mathfrak V}}
\def\C{{\mathfrak C}}
\def\X{{\mathfrak X}}
\def\Y{{\mathfrak Y}}
\def\R{{\bf R}}
\def\D{\widehat{D}}
\def\A{\widehat{A}}
\def\G{\hat{G}}
\def\T{{\cal T}}
\def\B{\widehat{B}}
\def\BC{\widehat{B_C}}
\def\BCC{\widehat{B_\C}}
\def\restr{\restriction}
\def\Aut{{\rm Aut\,}}
\def\Im{{\rm Im\,}}
\def\ker{{\rm ker\,}}
\def\inf{{\rm inf\,}}
\def\sup{{\rm sup\,}}
\def\Br{{\rm Br\,}}
\def\Yphi{Y_{[\phi]}}
\def\Ypsi{Y_{[\psi]}}
\def\Xphi{X_{\tilde{\phi}}}
\def\Xpsi{X_{\tilde{\psi}}}
\def\a{\alpha}
\def\abar{\overline{\alpha}}
\def\aa{{\bf a}}
\def\to{\rightarrow}
\def\arr{\longrightarrow}
\def\arl{\longleftarrow}
\def\sigmaa{{\bf \Sigma_a}}

\def\End{{\rm End\,}}
\def\Ines{{\rm Ines\,}}
\def\Hom{{\rm Hom\,}}

\def\restr{\upharpoonright}
\def\Ext{{\rm Ext}\,}
\def\Hom{{\rm Hom}\,}
\def\End{{\rm End}\,}
\def\Aut{{\rm Aut}\,}
\def\ker{{\rm ker}\,}
\def\defe{{\rm def}\,}
\def\rk{{\rm rk}\,}
\def\crk{{\rm crk}\,}
\def\nuc{{\rm nuc}\,}
\def\Dom{{\rm Dom}\,}
\def\Im{{\rm Im}\,}
\def\Yphi{Y_{[\phi]}}
\def\Ypsi{Y_{[\psi]}}
\def\Xphi{X_{\tilde{\phi}}}
\def\Xpsi{X_{\tilde{\psi}}}
\def\a{\alpha}
\def\abar{\overline{\alpha}}
\def\aa{{\bf a}}
\def\ra{\rightarrow}
\def\arr{\longrightarrow}
\def\iff{\Longleftrightarrow}
\def\sigmaa{{\bf \Sigma_a}}
\def\mm{{\mathfrak m}}
\def\F{{\mathfrak F}}
\def\X{{\mathfrak X}}
\def\Diam{\diamondsuit}
\def\mapdown#1{\Big\downarrow\rlap{$\vcenter{\hbox{$\scriptstyle#1$}}$}}
\def\cRp{\rm{cRep}$_2R$ }
\def\cR{c$R_2$}

\title{{\sc ALMOST FREE SPLITTERS}
\footnotetext{This work is supported by the project 
No. G-0294-081.06/93 of the German-Israeli
Foundation for Scientific Research \& Development\\
AMS subject classification:\\ 
primary 13C05, 18E40, 18G05, 20K20, 20K35, 20K40 \\
secondary: 13D30, 18G25, 20K25, 20K30, 13C10 \\
Key words and phrases: self-splitting modules, criteria for freeness of 
modules\\
GbSh 682 in Shelah's list of publications}
} 

\author{ R\"udiger G\"obel and Saharon Shelah}
\date{}

\maketitle

\begin{abstract}
Let $R$ be a subring of the rationals. We want to investigate self
splitting $R$-modules $G$ that is $\Ext_R(G,G) = 0$ holds. 
For simplicity we will call such modules splitters, see \cite{Sch}. Also other 
names like stones are used, see a dictionary in Ringel's paper \cite{Ri1}.
Our investigation continues \cite{GS}. In \cite{GS} we answered an open problem by 
constructing a large class of splitters. Classical splitters are free modules and 
torsion-free, algebraically compact ones. In \cite{GS} we concentrated on 
splitters which are larger then the continuum and such that countable 
submodules are not necessarily free. The `opposite' case of $\aleph_1$-free splitters of 
cardinality less or equal to 
$\aleph_1$ was singled out because of basically different techniques.
This is the target of the present paper. 
If the splitter is countable, then it must be free over some subring of the
rationals by Hausen \cite{Ha}. 
Contrary to the results in \cite{GS} and in accordance to \cite{Ha}
we can show that all $\aleph_1$-free splitters of cardinality $\aleph_1$
are free indeed. 
\end{abstract}

\section{Introduction}

Throughout this paper $R$ will denote a subring of the rationals $\Q$ 
and we will consider
$R$-modules in order to find out when they are splitters.
`Splitters' were introduced in Schultz \cite{Sch}. They also come up
under different names as mentioned in the abstract. 
\begin{definition} An $R$-module $G$ is a splitter if and only if 
$\Ext_R (G,G) = 0$ or 
equivalently if
$\Ext_{\Z} (G,G) = 0$ which is the case if and only if any 
$R$-module sequence
$$ 0 \arr G \overset{\beta}{\arr} X \overset{\alpha}{\arr} G \ra 0$$ 
splits.
\end{definition}

\medskip
\noindent
A short exact sequence
$$ 0 \arr B \overset{\beta}{\arr} C \overset{\alpha}{\arr} A \ra 0$$
represents $0$ in $\Ext(A,B)$ if and only if there is a splitting map
$\gamma: A \arr C$ such that $\gamma \alpha = id_A$. Here maps are 
acting on the right.

\medskip
\noindent
Recall an easy basic observation, see \cite{Fu}:

\medskip
\noindent
{\it If $\Ext (A,B) = 0, A' \subseteq A$ and $B' \subseteq B$, then 
$\Ext (A', B/B') = 0$ as well.} 

\medskip
\noindent
The first result showing freeness of splitters is much older then the notion of splitters and 
is due to Hausen \cite{Ha}. It says that any countable, torsion-free abelian group is a 
splitter if and only if 
it is free over its nucleus. The nucleus is the largest subring $R$ of $\Q$ which makes
the abelian group canonically into an $R$-module. 
More precisely

\begin{definition} The nucleus $R$ of a torsion-free abelian group
$G \neq 0$ is the subring $R$ of $\Q$ generated by all $\frac{1}{p}$
($p$ any prime) for which $G$ is $p$-divisible, i.e. $p G = G$.
\end{definition}

\medskip
\noindent
The fixed ring $R$ mentioned at the beginning will be the nucleus $R = \nuc G$ of
the associated abelian group $G$.

\medskip
\noindent
The following result reduces the study of splitters among abelian
groups to those which are torsion-free and reduced modules over their nuclei.

\begin{theorem}[\cite{Sch}]
Let $G$ be any abelian
group and $G = D \oplus C$ a decomposition of $G$ into the maximal divisible
subgroup $D$ and a reduced complement $C$. Then the following conditions are
equivalent.

\begin{enumerate}
\item $G$ is a splitter.
\item $\begin{cases} 
\mbox{(a)}\;\ D\;\ \mbox{is torsion (possibly 0),}\;\ C \;\
\mbox{is a torsion-free (reduced) splitter with}\\
\qquad p C = C \;\ \mbox{for all}\;\
p-\mbox{primary components}\;\ D_p \neq 0 \;\ \mbox{of}\;\ D.\\
\mbox{(b)}\;\ D \;\ \mbox{is not torsion and}\;\ C \;\ \mbox{is cotorsion.}
\end{cases}$
\end{enumerate}
\end{theorem}

Many splitters are constructed in
\cite{GS}, in fact we are also able to prescribe their endomorphism rings. This shows that 
uncountable splitters are not classifiable in any reasonable way, a result very much 
in contrast to classical well-known (uncountable) splitters which are the torsion-free 
algebraically compact (or cotorsion) groups.

\bigskip

The classical splitters come up naturally among many others when considering
Salce's work \cite{Sa} on cotorsion theories:
A  cotorsion theory is a pair of classes of $R$-modules $(\F, \C)$ 
which are maximal, closed under extensions such that the torsion-free class 
$\F$ is closed 
under subgroups and the cotorsion class $\C$ is closed under 
epimorphic images  
and $\Ext(F,C) = 0$ for all $F \in \F$ and $C \in \C$.
The elements in $\F \cap \C$ are splitters and
in case of Harrison's classical cotorsion theory these are the torsion-free, 
algebraically compact groups. For the trivial cotorsion theory these are 
 free $R$-modules.

Hausen's \cite{Ha} theorem mentioned above can be slightly extended
without much effort, see \cite{GS}.
\begin{theorem} \label{hausen}
If $R = \nuc G$ is the nucleus of the torsion-free
group $G$ and $G$ is a splitter of cardinality $< 2^{\aleph_0}$, then $G$ is
an $\aleph_1$-free $R$-module.
\end{theorem}
Recall that $G$ is an $\aleph_1$-free $R$-module if any countably generated
$R$-submodule is free. 
The algebraic key tool of this paper can be found in
Section 2. We consider torsion-free $R$-modules $M$ of finite rank which are 
minimal in rank and non-free. They are (by definition)  $n$-free-by-1 
$R$-modules if $\rk M = n+1$; the name is self explaining: 
They are pure extensions of a free $R$-module of rank $n$ by an $R$-module of
rank $1$. Similar to simply presented groups,
$n$-free-by-1 groups are easy represented by free generators and
relations. 
Using these minimal $R$-modules we will show the following
\begin{mtheorem} Any $\aleph_1$-free splitter of cardinality $\aleph_1$ 
is free over its nucleus.
\end{mtheorem}

The proof will depend on the existence of particular chains of 
$\aleph_1$-free $R$-module of cardinality $\aleph_1$ which we use to divide
$\aleph_1$-free $R$-modules of cardinality $\aleph_1$ into three 
types (I,II,III).
This may be interesting independently and we would like to draw attention 
to Section 3. In Sections 4-7 we use our knowledge about these chains to 
show freeness of splitter. The proof is divided into two main cases depending
on the continuum hypothesis CH (Section 5) and its negation (Section 4).
In the appendix Section 8 we present a proof of the main result of Section 5 
under the weaker set theoretic assumption WCH  $2^{\aleph_0} < 2^{\aleph_1}$, 
a weak form of CH 
which will be interesting (only) for splitters of cardinality $ > \aleph_1$.
The results in Section 6 and 7 on splitters of type II and III do not use
the case distinction by additional axioms of set theory.

\section{Solving Linear Equations}

Let $R$ be a subring of $\Q$. Then $R$-modules of minimal finite rank
which are not
free will lead to particular infinite systems of linear equations. 
Consider the Baer-Specker $R$-module $R^\omega$ of all $R$-valued functions
$f: \omega \rightarrow R$ on $\omega$, also denoted by 
$f = (f_m)_{m \in \omega}$.

\begin{lemma} \label{nosol}
Let $p = (p_m)_{m \in \omega}, \  
k_i = (k_{im})_{m \in \omega} \in R^\omega \ \ (i < n)$
where each $p_m$ is not a unit of $R$. Then we can find a sequence 
$s = (s_m)_{m \in \omega} \in R^\omega$ such that the following system of
equations $({\hat s})$ has no solution 
$\bar x = \ (x_0, \ldots, x_n) \in R^{n+1}, \ y_m \in R $ with 
$$y_0 = x_n , \\\ y_{m+1} p_m = y_m + \sum\limits_{i<n} x_i \; k_{im} + s_m \;\;\;
(m \in \omega). \ \ \ \ \;\;\; (\hat{s})$$
\end{lemma}

\noindent
{\bf Proof:} \; We will use Cantor's argument which shows that there are more real
numbers than rationals. First we enumerate all elements in $R^{n+1}$ as
$$ W = W = \{\overline x^m = (x^m_0, \ldots, x^m_n) \; : \; m \in \omega\}$$
and construct $s \in R^\omega$ inductively.

\medskip

It is interesting to note that the set of bad elements 
$$B = \{s \in R^\omega : \exists \bar{x} \in R^{n+1}, 
\ (y_m)_{m \in \omega} \in R^\omega \mbox { solving } \; 
(\hat{s})\} \subseteq R^\omega$$
is a submodule of $R^\omega$ but $|B|$ is uncountable in many cases. Hence 
enumerating $B$ would not help.

\medskip

Suppose $s_0, \ldots, s_{m-1} \in R$ are chosen and we must find $s_m$. We
calculate $y_0, \ldots, y_m$ from $s_0, \ldots, s_{m-1}$ and 
$y_0 = x^m_n, x^m_1, \ldots, x^m_{n-1}$ and equation $(\hat{s})$ up to $m-1$. The
values are uniquely defined by torsion-freeness and in particular\\
\begin{eqnarray} \label{sum}
z=y_m + \sum\limits_{i<n} x^m_i \; k_{im} \in R
\end{eqnarray}
is uniquely defined. Recall that $p_m$ is not a unit and either $p_m$ does
not divide $z$, then we set $s_m = 0$ or we can choose some
$s_m \in R \setminus \{ 0\}$
and $p_m$ does not divide $z + s_m$. In any case
\begin{eqnarray} \label{sum1}
 y \; p_m = z + s_m \mbox{ has no solution in } R
\end{eqnarray}
and $s_m$ is defined. 

Suppose that $(\hat{s})$ has a 
solution $\bar{x} \in R^{n+1}$, then
$\bar{x} = \overline x^m $ for some $m$ by our enumeration. We
calculate $y_{m+1}$ from $(\hat{s})$ substituting $\bar{x}
$, hence
$$y_{m+1} \; p_m = y_m + \sum\limits_{i < m} x_i^m \; k_{im} + s_m = z + s_m$$
is solvable by (\ref{sum}), which contradicts (\ref{sum1}). 
$\hfill \square$\\

\medskip

If $G' \subseteq G$ is a pure $R$-submodule of some $R$-module $G$ which is
of finite rank, not a free $R$-module such that all pure $R$-submodules of
$G'$ of smaller rank are free, then we will say that $G'$ is {\bf minimal
non-free}. Such minimal non-free modules are ``simply presented'' in the
sense that there are $x_i, y_m \in G'\ (i < n, m \in \omega)$
such that
\begin{eqnarray} \label{gr}
G' = \langle B, y_m R : m \in \omega \rangle \mbox{ with } 
B = \bigoplus\limits_{i < n} x_i R 
\end{eqnarray}
and the {\bf only relations} 
\begin{eqnarray} \label{rr}
y_{m+1} p_m = y_m + \sum\limits_{i < n} x_i \; k_{im}
\quad\quad (m \in \omega)
\end{eqnarray}
and coefficients $p_m,\; k_{im} \in R$. The submodule $B$ is pure in $G'$. 
If $G$ is not $\aleph_1$-free, then the existence of minimal non-free submodules
is immediate by Pontryagin's theorem.
Non-freeness of $G'$ implies that the Baer type of
$$G'/B = T \subseteq \Q$$
is strictly greater than the type of $R$, see 
Fuchs \cite[Vol 2, pp 107-112]{Fu}. 

In more details we have that $B$
is a pure submodule of $G'$, hence $G'/B$ is torsion-free of rank $1$ and
since $G'$ is not a free $R$-module, $G'/B = T$ cannot be isomorphic to $R$.
If $\varphi : G' \rightarrow \Q$ is the canonical homomorphism
taking $B$ to $0$ and $y_0$ to $1 \in \Q$, then $\Im \varphi \subseteq \Q$ 
represents the type of $\langle y_0 + B \rangle_* = G'/B$.
There are $ p_m \in \N$, not units in $R$ such that 
$T = \bigcup \limits_{m \in \omega} q^{-1}_m \Z \subseteq \Q$ 
 and $q_m = \prod_{i < m} p_i$. 
In order to derive the crucial 
equations as in the above definition we choose preimages 
$y_m \in G'$ of $q_m^{-1}$ such that
$$y_0 \varphi =q_0 = 1 \  \mbox{ and } y_m \varphi = q^{-1}_m 
\qquad \qquad (m \in \omega).$$
Using $q_{m+1} = q_m p_m$ we find elements 
$k_{i m} \in R \;\; (i < n),\;\; g_m \in G'$ 
such that (\ref{gr}) and (\ref{rr}) holds.
We will constantly use the representations (\ref{gr}) and (\ref{rr}) which are 
basic for the following

\begin{proposition}\label{embed}
Let $G_\alpha \subseteq G_{\alpha+1}$  be a countable free 
resolution of $G'$ as in (\ref{gr}) and let the relations (\ref{rr}) be expressed 
in $G_{\alpha +1}$ by
$$ y''_{m+1} p_m = y''_m + \sum\limits_{i < n} x''_i \; k_{im} + g_m$$
for some $g_m \in G_\alpha$; let $z_m$ $(m \in \omega)$ be non-trivial 
elements of an $\aleph_1$-free $R$-module $H^0$ of cardinality $\aleph_1$ 
and 
$$0 \arr H^0 \arr H_\alpha \overset{h}\arr G_\alpha \arr 0$$
be a short exact sequence.
Then we can find an $R$-module 
$$H' = \langle H_\alpha \oplus B', y'_{ m} : m \in \omega \rangle$$  
with $B'= \bigoplus\limits_{i<n} x'_iR$, $\overline g_m h = g_m$ and the 
only relations in $H'$
$$y'_{m+1} p_m = y'_m + \sum\limits_{i<n} x'_i k_{im} + z_m + \overline g_m.
\ \ \ \ (m \in \omega)$$
The map $h$ extends to $h'$ by $x'_ih' = x''_i,\ y'_mh' = y''_m$ such that  
the new diagram with vertical maps inclusions commutes:

\medskip\hfil\begin{tabular}{rcccccl}
$0 \arr$ & $H^0$         & $\arr$ & $H_\alpha$ & $\overset{h} 
\arr $ & $G_\alpha$ & $\arr 0$\\
         &$\downarrow$ && $\downarrow$ &                  & $\downarrow$ & \\
$0 \arr$ & $H^0$ & $\arr$ & $H'$ & $\overset{h'}
\arr$ & $G_{\alpha +1}$ & $\arr 0$\\
\end{tabular}\medskip\\
\end{proposition}

\medskip
{\bf Proof: } Let
$$F_{\alpha+1} = H_\alpha \oplus \bigoplus\limits_{i < n} \bar x_i R 
\oplus\bigoplus\limits_{m \in \omega} \bar y_m R $$
and define
$$N_{\alpha +1} = \langle (\bar y_{m+1} p_m - \bar y_m - \sum\limits_{i < n}
\bar x_i k_{im} - z_m - \bar g_m)R : m \in \omega \rangle.$$
Hence $H'= F_{\alpha + 1} / N_{\alpha + 1}$ and let
$$x'_i = \bar x_i + N_{\alpha + 1}, \; 
y'_m = \bar y_m + N_{\alpha + 1} \mbox{ and } x' = \bar x +N_{\alpha +1}.$$
First we see that\\

$(a)$ $\begin{cases}
x \arr x' \ (x \in H_\alpha) \mbox{ defines an embedding } H_\alpha \arr H'\\
\mbox{ and then we identify } H_\alpha \mbox{ with its image in } H'.
\end{cases}$\\

It remains to show that $H_ \alpha \cap N_{\alpha + 1} = 0$ viewed
in $F_{\alpha + 1}$. If $x \in H_\alpha \cap N_{\alpha + 1}$, then
there are $k_m \in R$ for $m \leq l$ and some $l \in \omega$ such that
$$\sum^l\limits_{m = 0} (\bar y_{m+1} p_m - \bar y_m - \sum\limits_{i < n} 
\bar x_i k_{i m} - z_m - \bar g_m) k_m = x \in H_\alpha.$$
We get $x = - \sum^l\limits_{m=0} ( z_m + \bar g_m) k_m$
and $\sum^l\limits_{m=0} (\bar y_{m+1} p_m - \bar y_m - \sum\limits_{i<n} 
\bar x_i k_{im}) k_m = 0.$
The coefficient of $\bar y_{l+1}$ is $p_l k_l = 0$, hence $k_l = 0$ and going 
down we get $k_m = 0$ for all $m \leq l$, hence $x = 0$ and (a) holds. 
Due to $N_{\alpha+1}$ we have the useful system of equations in $H'$.

$(b) \ \ \  y'_{m+1} p_m = y'_m + \sum\limits_{i<n} 
x'_i k_{im} + z_m + \bar g_m \mbox{ with } 
\bar g_m \in H_\alpha. $

In view of (a) we also have\\

$(c) \ \ \ H'= \langle H_{\alpha}\oplus B', y'_m :  m \in \omega \rangle 
\mbox{ with } B' = \bigoplus\limits_{i<n} x'_i R \subseteq H'.$

Next we claim that\\

$(d)$ $\begin{cases}
\mbox{ If } h'\upharpoonright H_\alpha = h, x'_i h'= x_i 
\mbox{ and } y'_m h'= y_m \ (i<n, m \in \omega), 
\mbox{ then }\\
h': H' \arr G_{\alpha+1} \mbox{ is a well-defined homomorphism with  }\\
\end{cases}$

\bigskip
$(e)$ $\ker h' = H^0$ and $\Im h' = G_{\alpha+1}.$\\

As $h'$ is defined on non-free generators, we must check that the relations 
between them are preserved, when passing to the proposed image. The relations 
are given by $N_{\alpha+1}$ or equivalently by $(b)$. Using the
definition (d) we see that the relations $(b)$ are mapped summand-wise 
under $h'$ as follows.

\medskip\hfil\begin{tabular}{rcccccl}
\quad & $y'_{m+1}p_m $ & $+$ & $y'_m$ $+$ & $\sum\limits_{i<n}
x'_i k_{im} $ $+$ & $\ \ z_m $ $+$ & $\bar g_m$\\
& $\downarrow$ & & $\downarrow$ & $\downarrow$ & $\downarrow$ & $\downarrow$\\
\quad & $y_{m+1}p_m $ & ? & $y_m$ $+$ & $\sum\limits_{i<n}x_i k_{im}$ 
$+$ & \;\, $0$ \, $+$ & $  g_m$\\
\end{tabular}\medskip

\noindent
and inspection of (\ref{rr}) and the relations in $G_{\alpha+1}$ shows that ? is an
equality sign ``=''. Hence $h'$ is well-defined. Notice that 
$H_\alpha h' = H_\alpha h = G_\alpha$, therefore $h'$ induces a
homomorphism
$$H'/H_\alpha \arr G_{\alpha + 1} / G_\alpha$$
and the last argument and $g_m \in G_\alpha$ show that this is an
isomorphism. Hence passing from $h$ to the extended map $h'$ the kernel 
cannot grow, we have $H^0 = \ker h' = \ker h $ and $\Im h'= G_{\alpha + 1}$ 
is obvious, so (d) and (e) and the proposition are shown.

\bigskip

\section{The Main Reduction Lemma - Types I, II and III}

The Chase radical $\nu G$ of a torsion-free $R$-module $G$ is the 
characteristic submodule
$$\nu G = \bigcap \{ U \subseteq G, \ G/U \mbox{ is } \aleph_1\mbox{-free } 
\}.$$
Since $G/\nu G$ is also $\aleph_1$-free, the Chase radical is the 
smallest submodule with $\aleph_1$-free quotient. 
If $U$ is a submodule of $G$ we write 
$$\nu_UG= G' \mbox{ for the Chase radical of } G \mbox{ over } U 
\mbox{ which is defined by } \nu (G/U) = G'/U. $$
Given any $\aleph_1$-free $R$-module $G$ of cardinality 
$|G| = \aleph_1$, we fix an $\aleph_1$-filtration 
$$G = \bigcup\limits_{\alpha < \omega_1} G^0_\alpha$$
which is an ascending, continuous chain of countable, free 
and pure $R$-submodule $G^0_\alpha$  of $G$ with $G^0_0 =0$.

We want to find a new ascending, continuous chain of pure $R$-submodules
$G_\alpha$ (not necessarily countable) such that 
$G = \bigcup\limits_{\alpha < \omega_1} G_\alpha$. However we do require that
\begin{eqnarray} \label{chase}
G/G_\alpha \mbox{ is } \aleph_1\mbox{-free if } \alpha \mbox{ is not a limit 
ordinal. } 
\end{eqnarray}
We will use the new chain to divide 
$\aleph_1$-free $R$-module of cardinality $ \aleph_1$ into three types. This
distinction helps to show that $\aleph_1$-free splitters of 
cardinality $\aleph_1$ are free.

Suppose $G_\beta \subseteq G$ is constructed for all $\beta < \alpha$. 
Next we want to define $G_\alpha$. 
If $\alpha$ is a limit ordinal, then 
$$G_\alpha = \bigcup\limits_{\beta <  \alpha} G_\beta.$$
Hence we may assume that $\alpha = \beta + 1$ and we must define
$G_\alpha =G_{\beta +1}$.
In order to ensure $G = \bigcup\limits_{\alpha \in \omega_1} G_\alpha$ 
we let
\begin{eqnarray}
G_{\alpha 0} = (G_\beta + G^0_\alpha)_* \subseteq G
\end{eqnarray}
the pure $R$-submodule generated by $G_\beta + G^0_\alpha$.
In any case we want to ensure that (\ref{chase}) holds, hence
$\nu_{G_{\alpha 0}}G \subseteq G_\alpha.$ Therefore we construct an 
ascending, continuous chain of pure $R$-submodules
\begin{eqnarray} 
\{ G_{\alpha j} : j < \omega_1 \} \mbox{ with } G_{\alpha j+1}/G_{\alpha j} 
= 0 \mbox{ or minimal non-free for each } 0 < j < \omega_1 
\end{eqnarray}
such that  $G_\alpha = \bigcup\limits_{j \in \omega_1} G_{\alpha j}$.
Suppose that $G_{\alpha i}$ is defined for all $i < j < \omega_1$.
If $j$ is a limit ordinal we take 
$G_{\alpha j} = \bigcup\limits_{i < j} G_{\alpha i}$ and if $j = i+1$ we 
distinguish two cases:
\begin{eqnarray} \label{count}
\mbox{If } G/G_{\alpha i} \mbox{ is } \aleph_1\mbox{-free, then }
G_{\alpha j} = G_{\alpha i} \mbox{ hence } G_\alpha = G_{\alpha i}  
\mbox{ and } |G_\alpha/G_{\alpha -1}| = \aleph_0.
\end{eqnarray}
Otherwise $G/G_{\alpha i}$ is not $\aleph_1$-free, and by Pontryagin's
theorem we can find a finite rank minimal non-free pure $R$-submodule 
$M/G_{\alpha i}$ of $G / G_{\alpha i}$ . Since 
$G = \bigcup\limits_{i \in \omega_1} G^0_i$ and 
$\emptyset \neq M \setminus G_{\alpha i} \leq G$, there is also a least 
ordinal $\gamma = \gamma (M) = \gamma (M/G_{\alpha i}) < \omega_1$, such that
\begin{eqnarray}
(M \setminus G_{\alpha i}) \cap (G^0_{\gamma + 1} \setminus 
G^0_\gamma) \neq \; \emptyset.
\end{eqnarray}

Among the candidates $M$ we choose one with the smallest $\gamma (M)$ 
and take it for $M = G_{\alpha i+1}$. This completes the construction of the 
$G_{\alpha i} \ 's$. Notice that either the construction of $G_\alpha$
stops as in case (\ref{count})  or we arrive at the second possibility:
\begin{eqnarray} \label{uncount}
G_{\alpha i+1}/G_{\alpha i} \mbox{ is minimal non-free for each } 
i < \omega_1  \mbox{ and } |G_\alpha/G_{\alpha-1}| = \aleph_1.
\end{eqnarray}

It remains to show that in case (\ref{uncount}) the following holds.
\begin{eqnarray} \label{free}
\nu_{G_{\alpha 0}}G = G_\alpha  \mbox{ or equivalently }
G/G_\alpha \mbox{ is } \aleph_1 \mbox{-free.}
\end{eqnarray}
Suppose that $G/G_\alpha$ is not $\aleph_1$-free and let $X$ be a
non-free submodule of minimal finite rank in $G/G_\alpha$ which exists by 
Pontryagin's theorem. 
Representing $X$ in $G$ we have
$$G''= \langle x_i, y_m, G_\alpha : i < n, m \in \omega\rangle_* 
\mbox{ with } G''/ G_\alpha = X,$$ 
see also G\"obel, Shelah \cite{GS}. 
There are elements $g_m \in G_\alpha \ (m \in \omega)$ such that 
$$y_{m+1} p_m = y_m + \sum\limits_{i < n} x_i k_{im} + g_m$$
for some $p_m, k_{im} \in R$ ($p_m$ not units of $R$). We take 
$$G' = \langle x_i, y_m, g_m : i < n, m \in \omega \rangle_* \subseteq G,$$ 
hence $X = G' + G_\alpha / G_\alpha$ was our starting point. Since $G'$ is 
obviously countable, there is a $\gamma^* \in \omega_1$ with $G' \subseteq  
G^0_{\gamma^*}$. 
If $\{ G_{\alpha j} : j < \omega_1 \}$ is the chain constructed above,
we also find $i \in {\omega_1}$ with $g_m \in G_{\alpha i}$ 
for all ${m \in \omega}$.
If $i \leq j \in \omega_1$, then $G' + G_{\alpha j} / G_{\alpha j}$ is an
epimorphic image of $X$, hence minimal non-free or $0$.
The second case leads to the immediate contradiction:
$$G' \subseteq G_{\alpha j} \subseteq G_\alpha \;\;\mbox{but}\;\;X \neq 0.$$
Hence $G' + G_{\alpha j} / G_{\alpha j} \neq 0$ was a candidate for 
constructing $G_{\alpha j+1}$ for any $i \leq j \in \omega_1$. 
Has it been used? 
We must compare the $\gamma$-invariant 
$\gamma (G' + G_{\alpha j}/G_{\alpha j})$ with the various 
$\gamma (G_{\alpha j+1}/G_{\alpha j})$. From $G' \subseteq G^0_{\gamma^*}$ 
we see that there is 
$\gamma^j < \gamma^*$ such that 
$$(G' + G_{\alpha j} \setminus G_{\alpha j}) \cap (G^0_{\gamma^j+1} 
\setminus G^0_{\gamma^j} \neq \emptyset.$$
By minimality of $\gamma_j = : \gamma (G_{\alpha j+1} / G_{\alpha j})$ 
we must have $\gamma_j \leq \gamma^j < \gamma^*$ and 
$$(G_{\alpha j+1} \setminus  G_{\alpha j}) \cap (G^0_{\gamma_j+1} \setminus 
G^0_{\gamma_j}) \neq \emptyset $$
and $(G_{\alpha j} \cap G^0_{\gamma^*}) \quad (j \in \omega_1)$
is a {\bf strictly} increasing chain of length $\omega_1$ of the 
countable module $G^0_{\gamma^*}$, which is impossible.
Hence $G/G_\alpha$ is $\aleph_1$-free and (\ref{free}) is shown. 
We have a useful additional property of the constructed chain which reflects 
(\ref{free}).
\begin{corollary} \label{nu}
If $0 \neq \alpha \in \omega_1 $ is not a limit ordinal, then 
$G_\alpha = \nu_{G^0_\alpha}G$.
\end{corollary}
{\bf Proof:} 
We concentrate on the case (\ref{uncount}) and only note that the case 
(\ref{count}) is similar.

Recall from (\ref{free}), that $G/G_\alpha$ is $\aleph_1$-free,
hence the claim of the corollary is equivalent to say that any submodule $X$
of $G_\alpha$ must be  $ G_\alpha$  if only $G^0_\alpha \subseteq X$ with 
$G_\alpha /X \ \aleph_1$-free.
 
Let $\alpha > 0$ and suppose 
$G^0_\alpha \subseteq X \subseteq G_\alpha$ and $0 \neq G_\alpha / X$ is 
$\aleph_1$-free. First we claim that 
$$G_\beta \subseteq X \mbox{ for all } \beta < \alpha.$$ 
If this is not the case, then let $\beta < \alpha$ be minimal with $G_\beta 
\not\subseteq X$. Recall that $\beta $ can not be a limit
ordinal and we can write $\beta = \gamma + 1$ for some $\gamma < \beta$.
We have $G_\beta = \bigcup\limits_{j \in \omega_1} \; G_{\beta j}$, hence
$$i_\beta = \mbox{min}\; \{j \in \omega_1 : G_{\beta j} \not\subseteq X \} 
\in\omega_1$$
exists. If $i_\beta = 0$, then 
$G^0_\beta \subseteq G^0_\alpha \subseteq X$ from $\alpha > \beta$ and 
$G_{\beta 0} \not\subseteq X$.   
We get $G_{\beta 0} = \langle G_\gamma ,G^0_\beta \rangle_* 
\not\subseteq X$ and $G^0_\beta \subseteq X$ requires
$G_\gamma \not\subseteq X$, contradicting minimality of $\beta$. 
Hence $i_\beta > 0$ and $i_\beta = j + 1$. 
We have $G_{\beta i_\beta} \not\subseteq X$ and 
$G_{\beta j} \subseteq X$ from $j < i_\beta$ and minimality of 
$i_\beta$. However $ G_{\beta j+1}/G_{\beta j}$ is minimal non-free, and
$0 \neq G_{\beta j+1} + X/X \subseteq G/X$
is an epimorphic image, hence non-free as well. Therefore $G/X$ is not
$\aleph_1$-free, a contradiction showing our first 
claim.\\
From the first claim we derive $\bigcup\limits_{\beta < \alpha} 
G_\beta\subseteq X$. Now there must be a minimal
$$i_\alpha = \;\mbox{min}\; \{j : G_{\alpha j} \not\subseteq X\} \in 
\omega_1,$$
which cannot be a limit ordinal, and again $i_\alpha > 0$, 
hence $i_\alpha = j+1$. We find 
$G_{\alpha j} \subseteq X, G_{\alpha j+1} \not\subseteq X$ and 
$G_\alpha / X$ cannot be $\aleph_1$-free, a final contradiction.
$\hfill \square$\\

We now distinguish cases for $G$ depending on the existence of
particular filtrations. Let 
$G= \bigcup\limits_{\alpha < \omega_1} G_\alpha$ be the filtration 
constructed from the $\aleph_1$-filtration 
$G= \bigcup\limits_{\alpha < \omega_1} G^0_\alpha$. 

If there is an ordinal 
$\beta < \omega_1$ (which we assume to be minimal) such that $G = G_\beta$, 
then let $C = G^0_\beta$ which is a countable, free and pure $R$-submodule 
of $G$. From Corollary \ref{nu} we see that $\nu_CG = G$. Hence, beginning 
with $C$ we get a new $\aleph_1$-filtration (we use the same notation)
$\{G_\alpha : \alpha \in \omega_1\}$ of countable, pure and free 
$R$-submodules of $G$ such that $G_0 = C$ and 
{\bf each} $G_{\alpha + 1} / G_\alpha \ (\alpha > 0)$ is minimal non-free. 
In this case we say that $G$ and the filtration are of type I.

In the opposite case the chain only terminates at the limit 
ordinal $\omega_1$ i.e. $G_\beta \neq G$ for all $\beta < \omega_1$.
We have a proper filtration $G = \bigcup\limits_{\alpha < 
\omega_1}G_\alpha$ such that Corollary \ref{nu} holds. 
If for each $\alpha \in \omega_1$ for some $i < \omega_1$ case (\ref{count})
occurs, then the constructed chain $\{G_\alpha : \alpha \in \omega_1\}$ is 
an $\aleph_1$-filtration of countable, pure and free $R$-submodules with the 
properties of Corollary \ref{nu} and (\ref{chase}). We say that the chain 
and $G$ are of type II respectively.

If $G$ is not of type I or of type II we say that
$G$ is of type III. In this case, there is a first
$\alpha \in \omega_1$ such that $G_{\alpha+1}/G_\alpha$ is uncountable. We
may assume that $\alpha = 0$. With the new enumeration we see that
the following holds for type III:

$(III) \ \begin{cases}

G= \bigcup\limits_{\alpha < \omega_1} G_\alpha, \; G_0 = 0, 
|G_1| = \aleph_1 \mbox{ and (3.1) holds, }\\
G_1 = \bigcup\limits_{j \in \omega_1} G_{0 j} \mbox{ is an } 
\aleph_1 \mbox{-filtration of pure submodules of } G_1\\
\mbox{ with \bf{each} } G_{0 j+1} / G_{0 j} \mbox{ minimal non-free}.
\end{cases}$

\bigskip
We have a
\noindent
\begin{reduction-lemma} \label{type}
Any $\aleph_1$-free module $G$ of cardinality
$\aleph_1$ is either of type I, II or III.
\end{reduction-lemma}

\bigskip
\section{Splitters Of Cardinality $\aleph_1 < 2^{\aleph_0}$ Are Free}
In this section we do not need the classification of $\aleph_1$-free 
$R$-modules of cardinality $\aleph_1$ given in Lemma \ref{type}. 
Moreover, we note that $\aleph_1$-freeness of splitter of cardinality
$\aleph_1 < 2^{\aleph_0}$ follows by Theorem \ref{hausen}. In fact we
will present a uniform proof showing freeness of splitters up to cardinality
$\aleph_1 < 2^{\aleph_0}$ which extends Hausen's
result \cite{Ha} concerning countable splitters.
 We begin with a trivial observation

\begin{proposition}  \label{nuc}
Let $G = \bigcup\limits_{\alpha \in \omega_1} G_\alpha$
be an $\aleph_1$-filtration of pure and free $R$-submodules 
$G_\alpha$ of $G$. Then $\nuc ( G_\alpha) = R$ for all 
$\alpha \in \omega_1$.
\end{proposition}

\noindent
{\bf Proof:}  Choose any basic element $b \in G_\alpha$ for some
$\alpha \in \omega_1$. If $r \in \Q$ divides $b$ in $G$, 
then $r$ divides $b$ in $G_\alpha$ by purity, hence $r \in R$ from 
$bR \oplus C = G_\alpha$ and $\nuc G = R$.\\

\begin{corollary} \label{nuc1}
$(\aleph_1 < 2^{\aleph_0})$. If $G$ is a splitter of
cardinality $\leq \aleph_1$ and $\nuc G = R$, then there is an 
$\aleph_1$-filtration $G = \bigcup\limits_{\alpha \in \omega_1}G_\alpha$ of 
pure and free $R$-submodules $G_\alpha$ such that $\nuc G_\alpha = R$ 
for all $\alpha \in\omega_1$.  
\end{corollary}

\bigskip
\noindent
{\bf Proof:} From $\aleph_1 < 2^{\aleph_0}$ and G\"obel and Shelah 
\cite{GS}, see Theorem \ref{hausen}, 
follows that $G$ is an $\aleph_1$-free $R$-module and $G$ has 
an $\aleph_1$-filtration as in the hypothesis of Proposition \ref{nuc}. 
Hence Corollary \ref{nuc1} follows by the Proposition \ref{nuc}.
 $\hfill \square$

\begin{definition} \label{White}
Let $G$ be a torsion-free abelian group with
$\nuc G = R$ and $X$ an $R$-submodule of $G$. Then $X$ is contra-Whitehead in
$G$ if the following holds. \\
There are $z_m \in G,\ p_m, \ k_{im} \in R \ (i < n, \ m \in \omega)$ 
such that the system of equations 
$$Y_{m+1} p_m \equiv Y_m + \sum\limits_{i < n} X_i k_{im} + z_m 
\mbox{\it  \  mod \  X }
\qquad \ \ (m \in \omega)$$
has no solutions $y_m, a_i \in G$ (for $Y_m, X_i$ respectively) with 
$\bigoplus\limits_{i < n} (a_i + X) R$ free of rank $n$ and pure in 
$G/X$. Otherwise we call $X$ pro-Whitehead in $G$.
\end{definition}

For $X \subseteq G$ as in the definition let 
$ \mathfrak W $ be the set of all finite
sequences $\bar{a} = (a_0, a_1 \ldots, a_n)$ such that \\
(i) \hspace*{0,5cm} $a_i \in G \;\; (i \leq n)$\\
(ii)\hspace*{0,4cm} $\bigoplus\limits_{i < n} (a_i + X)\; R$ 
is pure in $G/X$.\\
(iii)\hspace*{0,4cm} $\left <(a_i + X)R : i \leq n \right>_*$ 
is not a free $R$-module in $G/X$.

\noindent
In particular $G'_{\bar{a}} = \bigoplus\limits_{i < n} a_i R\; \oplus X$ is
a pure submodule of 
$G_{\bar{a}} = \langle X, \ a_i R : i \leq n \rangle_*$ and
of $G$, the module $G_{\bar{a}}/X$ is an $n$-free-by-1 $R$-module.
From (\ref{rr}) we find $p_{\bar a m} \in \N$ not units in $R$ and elements 
$k_{\bar{a}i m} \in R \;\; (i < n),\;\; g_{\bar a m} \in G_{\bar{a}}$ 
such that
\begin{eqnarray} \label{equab}
y_{\bar a m+1} p_{\bar a m} = y_{\bar a m} + 
\sum\limits_{i < n} a_ik_{\bar a im} + g_{\bar{a}m} \qquad (m \in \omega)
\end{eqnarray}

The equations (\ref{equab}) are the basic systems of equations which decide
about $G$ to be a splitter. We will also consider an `inhomogeneous counter 
part' of (\ref{equab}) and choose a sequence 
$\bar z = (z_m : m \in \omega)$ of elements $z_m \in G$. 
The $\bar z$-inhomogeneous counter part of (\ref{equab}) is the system of 
equations
\begin{eqnarray} \label{inho}
Y_{m+1} p_{\bar a m} \equiv Y_m + \sum\limits_{i < n} X_ik_{\bar a im} + z_m 
\mbox{\it \ mod \ X }
\qquad \ \  (m \in \omega)
\end{eqnarray}

\noindent
According to the above definition we also say that $\bar a \in \mathfrak W$ 
is contra-Whitehead if (\ref{inho}) has no solutions $y_m \ (m \in \omega)$ 
in $G$ (hence in $G_{\bar a})$ for some $\bar{z}$ and $X_i = a_i$. 
Otherwise we say that $\bar a$ is pro-Whitehead.
If $G = \bigcup\limits_{\alpha \in \omega_1} G_\alpha$ is an 
$\aleph_1$-filtration of $G$, then we define $ {\mathfrak W}_\alpha$ for 
$X = G_\alpha$ 
and let $S = \{\alpha \in \omega_1$: there exists 
$\bar a \in {\mathfrak W}_\alpha$ 
contra-Whitehead $\}$.

\begin{proposition} \label{stat}
If $G = \bigcup\limits_{\alpha\in \omega_1} G_\alpha$ 
and $S$ as above is stationary in $\omega_1$, then $G$ is not a splitter. 
\end{proposition}

Before proving this proposition we simplify our notation. If $\alpha \in S$ 
we choose 
$$z^\alpha_m \in G,\ \bar a = (a^\alpha_0, a^\alpha_1, \ldots, a^\alpha_n),
\ p_{\bar a m} = p_{\alpha m}, \ g_{\bar a m} = g^\alpha_m, \
k_{\bar a im} = k_{\alpha i m}, \ y_{\bar a m} = y^\alpha_m$$
so that equations (\ref{equab}) and (\ref{inho}) become for $X = G_\alpha$ 
\begin{eqnarray} \label{equaa}
y^\alpha_{m+1} p_{\alpha m} = y^\alpha_m + \sum\limits_{i < n} 
a^\alpha_ik_{\alpha i m} + g^\alpha_m \quad (m \in \omega)
\end{eqnarray}
with $\bar{z}^\alpha$-inhomogeneous counter part
\begin{eqnarray} \label{inhoa} 
Y_{m+1} p_{\alpha m} \equiv Y_m + \sum\limits_{i < n} X_ik_{\alpha i m} + 
z^\alpha_m \mbox{\it \ mod \ } G_\alpha
\qquad (m \in \omega).
\end{eqnarray}
Hence (\ref{equaa}) is a system of equations with solutions 
$y^\alpha_m, \ a^\alpha_i, \ g^\alpha_m$ in $G_{\alpha + 1}$, 
while (\ref{inhoa}) with variables
$Y_m,\  X_i \ (m \in \omega,\ i < n)$ has no solutions in $G$,
as discussed in Definition \ref{White} for $X = G_\alpha$.
The set of limit ordinals is a cub, hence we may restrict $S$ to this cub
and assume that $S$ consists of limit ordinals only. 
If $\alpha \in S$ we also may assume that
$$G_{\alpha+1} = \langle G_\alpha,\ a^\alpha_i R : i < n \rangle_* 
=\langle G_\alpha, \ a^\alpha_i R,\  y^\alpha_m R : m \in \omega_1, 
\;\; i < n  \rangle.$$
We begin the

\noindent
{\bf Proof of the Proposition \ref{stat}:} 

\noindent
We will use the last remarks 
for constructing $h : H \arr G$ such that
$$ (*)\qquad  0 \arr H^0 \arr H \overset{h} \arr G \arr 0$$
does not split, hence $\Ext (G, H^0) \neq 0$. We will have 
$H^0 \cong G$, hence $\Ext (G, G) \neq 0$ and $G$ is not a splitter.
\\
Choose an isomorphism $\gamma : G \arr H^0$ which carries the 
$\aleph_1$-filtration 
$\{G_\alpha : \alpha \in \omega_1\}$ to 
$H^0 = \bigcup\limits_{\alpha \in \omega_1} H'_\alpha$ and 
$z^\alpha_m$ to ${z'}_{\alpha m}$. 
Inductively we want to define short exact sequences
$$(\beta) \qquad 0 \arr H^0 \overset{id} \arr H_\beta \overset{h_\beta} 
\arr G_\beta \arr 0 \qquad\;\; (\beta < \alpha)$$
which are increasing continuously. Let
$$(0) \qquad 0 \arr H^0 \overset{id} \arr H_0 \overset{h_0} \arr 0 \arr 
0$$
be defined for $H_0 = H^0$ with $h_0$ the zero-map and suppose $(\beta)$ 
is defined for all $\beta < \alpha < \omega_1$ with $\alpha$ a limit 
ordinal. We take unions and $(\alpha)$ is defined. If $\alpha \in \omega_1 
\setminus S$, we extend $(\alpha)$ trivially to get $(\alpha + 1)$ and if 
$\alpha \in S$ we must work for $(\alpha + 1)$:\\
We apply Proposition \ref{embed} to find 
$H_\alpha \subseteq H_{\alpha + 1}$
with 
$$H_{\alpha+1} = \langle H_\alpha, e_{\alpha m}, \ x_{\alpha m}\; :\; 
m \in \omega,\;\; i < n\rangle.$$
and relations
\begin{eqnarray} \label{Haa}
e_{\alpha m+1} \; p_{\alpha m} = e_{\alpha m} + \sum\limits_{i < n} 
x_{\alpha i} k_{\alpha i m} + y_{\alpha m} + z'_{\alpha m} 
\;\;\; \;(m \in \omega)
\end{eqnarray}
with $y_{\alpha m} h_\alpha = g^\alpha_m \in G_\alpha$. 
We want to extend the homomorphism $h_\alpha: H_\alpha \arr G_\alpha$ 
to $h_{\alpha + 1} : H_{\alpha + 1} \arr G_{\alpha+1}$, and set 
$e_{\alpha m} h_{\alpha + 1} = y^\alpha_m$ and $x_{\alpha i} h_{\alpha + 1} 
= a^\alpha_i$. By Proposition \ref{embed} the map 
$h_{\alpha + 1}$ is a well defined 
homomorphism. It is clearly surjective with kernel $H^0$. Hence $(\alpha + 
1)$ is well-defined for all $\alpha \in\omega_1$ and 
$h = \bigcup\limits_{\alpha \in \omega_1} h_\alpha$ shows $(*)$.\\ 
Finally we must show that $(*)$ does not split and suppose that 
$\sigma : G \arr H$ is a splitting map for $(*)$, hence 
$\sigma h = \mbox{id}_G$ and $H = H^0 \oplus \mbox{Im}\sigma$ and
$g^\alpha_m = y_{\alpha m} h_\alpha = y_{\alpha m} h $, so 
$(y_{\alpha m}-g^\alpha_m \sigma)\; h = 0$ implies 
$y_{\alpha m} - g^\alpha_m \sigma \in H^0$ for all $\alpha \in S$. The set 
$$C = 
\{ \alpha \in \omega_1 : \alpha \mbox{ a limit ordinal } y_{\alpha m} - 
g^\alpha_e \ \sigma \in H'_\alpha \}$$ 
- by a back-and-forth argument - is a cub 
and hence $S \cap C$ is stationary in $\omega_1$. We can find 
$\alpha \in C \cap S$ 
and consider the attached equations. In $G$ holds  (\ref{equaa})
$$y^\alpha_{m+1}\ p_{\alpha m} = y^\alpha_m + \sum\limits_{i < n} 
a^\alpha_ik_{\alpha i m} + g^\alpha_m$$
and $\sigma$ moves these equations to $H$:
$$(y^\alpha_{m+1}\; \sigma)\; p_{\alpha m} = 
(y^\alpha_m \; \sigma) + \sum\limits_{i < n}\; (a^\alpha_i \;\sigma)\; k_{\alpha i m} 
+ (g^\alpha_m\; \sigma)$$
which we subtract from (\ref{Haa}). 
Hence

$$(e^\alpha_{m+1} - y^\alpha_{m+1}\;\sigma)\; p_{\alpha m} = 
(e^\alpha_m - y^\alpha_m \;\sigma)
+ \sum\limits_{i<n} (x_{\alpha i} - a^\alpha_i\;\sigma)\; k_{\alpha i m} + 
(y_{\alpha m} - g^\alpha_m\;\sigma) + z'_{\alpha m}$$
Put 
$$f_{\alpha m} = e^\alpha_m - g^\alpha_m \sigma, \ 
v_{\alpha i} = x_{\alpha i}- a^\alpha_i \sigma, \ 
w_{\alpha m} = y_{\alpha m} - g^\alpha_m \;\sigma$$
and note that
$$f_{\alpha m} h = e_{\alpha m} h - g^\alpha_m \sigma h = 
g^\alpha_m - g^\alpha_m= 0$$ 
hence $f_{\alpha m} \in \ker h = H^0$. Similarly $w_{\alpha m}, \ 
v_{\alpha m} \in H^0$. The last equation turns into
$$f_{\alpha m+1} \ p_{\alpha m} = f_{\alpha m} + \sum\limits_{i<n} \ 
v_{\alpha i m} \ k_{\alpha i m} + w_{\alpha m} + z'_{\alpha m}
\ \ \ (m \in \omega)$$
which, as just seen, is a system of equations in $H^0$. 
From $\alpha \in C$ we have $w_{\alpha m} \in H'_\alpha$. 
The isomorphism $\gamma^{-1}$ moves the last equation back into $G$ and 
$w_{\alpha m} \gamma^{-1} \in G_\alpha$. Using 
$$f'_{\alpha m} = f_{\alpha m} \gamma^{-1},\ w'_{\alpha m} = w_{\alpha m} 
\gamma^{-1}, \ v'_{\alpha i} = v_{\alpha i} \gamma^{-1}$$ 
we derive 
$$f'_{\alpha m+1} p_{\alpha m} = f'_{\alpha m} +\sum\limits_{i<n} 
v'_{\alpha i} k_{\alpha i m} + w'_{\alpha m} + z_\alpha^m \ \ \ 
(m \in \omega)$$
with $w'_{\alpha m} \in G_\alpha$ and $z_\alpha^m$ as in 
(\ref{inhoa}), which is impossible in the case $\alpha \in S$ 
which is contra-Whitehead, where we have chosen $z_\alpha^m$ 
accordingly. $\hfill \square$

\begin{theorem} \label{afree}
Let $G$ be a splitter of cardinality  $< 2^{\aleph_0}$ with 
$\nuc G = R$. If $X$ is a pure, countable $R$-submodule of $G$ which is  
pro-Whitehead in $G$, then $G/X$ is an $\aleph_1$-free $R$-module.
\end{theorem}
\noindent
{\bf Proof:} First we assume that $\nuc(G/X) = R$ and  
suppose for contradiction that $G/X$ is not an $\aleph_1$-free $R$-module. 
By Pontryagin's theorem we can find an $R$-submodule $Y \subseteq G/X$ of 
finite rank which is not free. We may assume that $Y$ is of minimal 
rank. Hence
$$Y = \langle B, y_m R : m \in \omega \rangle,\;\;\; B = 
\bigoplus\limits_{i<n} \bar{x}_i R$$
with the only relations 
$$y_{m+1} p_m = y_m + \sum\limits_{i<n} \bar{x}_i k_{i m}\quad (m \in 
\omega)$$ 
as in Section 2 such that each $p_m \in R$ is 
not a unit of $R$ for $m \in \omega$. Choose
$x_i \in G$ such that $x_i + X = \bar{x}_i$ for each $i < n$. We can 
also choose a sequence of elements 
$$z_m \in G \mbox{ such that } z_m + X \mbox{ is not divisible by } 
p_{m-1} \mbox{ from nuc}(G/X) = R \ \ (m \in \omega).$$ 
If $\eta \in {}^\omega 2$, then let
$$\bar{z}^\eta = \langle \eta (e) z_e : e \in \omega \rangle =(z^\eta_e).$$ 
Recall that $X$ is pro-Whitehead in $G$, hence the systems of equations

\begin{eqnarray*} 
y^\eta_{m+1} p_m \equiv y^\eta_m + \sum\limits_{i<n} x^k_i k_{i m}+ z^\eta_m 
\mbox{\it \ mod \ X } \qquad (m \in \omega) \ \ \qquad \qquad  ( \eta )
\end{eqnarray*}

\noindent
has solutions $x^\eta_i$, $ y^\eta_m \in G$ for each $\eta \in {}^\omega 2$.
Note that 
$$| \{ \langle x^\eta_i : i < n \rangle  {}^\wedge \langle y^\eta_0
\rangle :\eta \in {}^\omega 2 \}| \leq |G| < 2^{\aleph_0}.$$ 
We can find $\eta \neq \nu \in {}^\omega 2$ such that 
$x^\eta_i = x^\nu_i$ for all $i < n$ and $y^\eta_0 = y^\nu_0$. 
From $\eta \neq \nu$ we find a branching point $j \in \omega$ such that 
$$ \eta (j) \neq \nu (j)  \mbox{ but } 
\eta \upharpoonright j = \nu\upharpoonright j.$$
We may assume 
$$\eta (j) = 1 \mbox{ and } \nu (j) = 0 $$
and put $w_m = y^\eta_m - y^\nu_m$. Subtracting the equations 
($\nu$) from ($\eta$) we get from 
$x^\eta_i - x^\nu_i = 0$ that
$$w_{m+1}\; p_m = w_m + (z^\eta_m - z^\nu_m) \mbox{\it \ mod \ X }$$
and $w_0 = y^\eta_0 - y^\nu_0 = 0$ as well. For $m \leq j$ we have 
$z^\eta_m - z^\nu_m = 0$ and $z^\eta_j - z^\nu_j = z_j$, hence $w_m = 0$ for 
$m < j$ by torsion-freeness and
$$w_j \;\;p_{j-1} = z_j \mbox{ \it \ mod \ X }$$
which contradicts our choice of $z_m$'s and $p_m$s.

If $\nuc(G/X) =\Q$, then $G/X$ is divisible, hence $X$ is dense and pure in
$G$, we have 
$X \subseteq_* G \subseteq_* \widehat X$, where $\widehat X$ is 
the $\Z$-adic completion of  $X$, and $X$ is a free $R$-modules of countable
rank. Hence 
$G/X \subseteq_*\widehat X/X \cong \bigoplus\limits_{2^{\aleph_0}} \Q$ and 
there are $2^{\aleph_0}$ independent elements in $\widehat X/X$.
Using these independent elements, we find systems of equations expressing
them as solutions $\Z$-adic limits, which must be solvable by pro-Whitehead. 
Hence $|G| = |\widehat X|= 2^{\aleph_0}> \aleph_1$, which is a contradition.

So we find $p_m\in R$ and $z_m \in G$  such that $p_{m-1}$ does not devide 
$z_m + X$ in $G$ for all $m \in \omega$. 
The above argument applies again for $n=0$ and leads to a 
contradiction.

\begin{corollary} \label{nonCH}
Any splitter of cardinality at most $\aleph_1 < 2^{\aleph_0}$
is free over its nucleus.
\end{corollary}
{\bf Proof:} Let 
$$G = \bigcup\limits_{\alpha \in \omega_1} G_\alpha$$  
be an $\aleph_1$-filtration of the splitter $G$. By Corollary \ref{nuc1}
we may assume that each $G_\alpha$ is a pure and free $R$-submodule of $G$ 
with $R = \nuc G$. If $S$ denotes the set 
$$\{ \alpha \in \omega_1 : G_\alpha \mbox{ is contra-Whitehead in } G \},$$ 
then $S$ is not stationary in $\omega_1$ by the last Proposition \ref{stat}. 
We may assume that all $G_\alpha$ are pro-Whitehead in $G$ and each 
$G_{\a+1} / G_\a$ is countable, hence free by Theorem \ref{afree}. 
We see that $G$ must be free as well. $\hfill \square$
\bigskip
\section{Splitters Of Type I Under CH}
In view of Section 4 we may assume CH to derive a theorem in ZFC showing 
freeness for $\aleph_1$-free splitters of cardinality $\aleph_1$ of type I.
The advantage of the set theoretical assumption is - compared with the 
proof based on the weak continuum hypothesis WCH in Section 8 - that the 
proof given here by no means is technical. 
Recall that $G$ is of type I if 
$G = \bigcup\limits_{\alpha \in \omega_1} G_\alpha$
for some $\aleph_1$-filtration $\{G_\alpha : \alpha \in \omega_1\}$ of pure
submodules $G_\alpha$ such that {\bf each } $G_{\alpha+1} / G_\alpha
\ (\alpha > 0)$ is a minimal non-free $R$-module. 
In this section we want to show the following
\begin{proposition} \label{typeI}
{\rm (ZFC + CH)} \ \ Modules of type I are not splitters. 
\end{proposition}

Combining Proposition \ref{typeI} and Corollary \ref{nonCH} we have 
can remove CH and have the 
immediate consequence which holds in ZFC.

\begin{corollary} 
Any $\aleph_1$-free splitter of type I (and cardinality $\aleph_1$) is free 
over its nucleus.
\end{corollary}

The proof of Proposition \ref{typeI} is based on an observation strongly 
related to type I concerning splitting maps. 
Then we want to prove a step lemma for applications of CH.
Finally we use CH to show  $\Ext (G,G) \neq 0$ in Theorem \ref{typeI}.
In Section 1 we noticed that if 
$$0 \arr B \overset{\beta}{\arr} C \overset{\alpha}{\arr}A \arr 0$$
is a short exact sequence, hence representing an element in $\Ext (A,B)$,
then this element is $0$ if and only if there is a splitting map $\gamma:
A \arr C$ such that $\gamma \alpha = id_A$. This simple fact is the key for 
the next two results.
\begin{observation}  \label{obs1}
Let 
$G = \bigcup\limits_{\alpha \in \omega_1}G_\alpha$ be a
filtration of type I. For ${\alpha \in \omega_1}$, let
$$0 \arr H^0 \arr H_\alpha \arr G_\alpha \arr 0 \quad\quad\quad (\alpha)$$
be a continuous, increasing chain of short exact sequences with union
$$0 \arr H^0 \arr H \arr G \arr 0 \quad\quad\quad (\omega_1)$$
and let $H^0 \cong G$ be $\aleph_1$-free. Then any splitting map of $(1)$ 
has at most one extension to a splitting map of $(\omega_1)$.
\end{observation}

\noindent
{\bf Proof:} \ We may assume that the splitting map 
$\sigma: G_1 \arr H_1$ of $(1)$ has two extensions 
$\sigma, \sigma': G \arr H$ which split. 
Since $\omega_1$ is a limit ordinal, there is some $\beta < \omega_1$ minimal 
with $(\sigma - \sigma') \upharpoonright G_{\beta} \neq 0$. 
Clearly $\beta $ is not a limit ordinal and 
$\sigma - \sigma'$ induces a non-trivial map
$\delta: G_{\beta}/G_{\beta-1} \arr H_\beta.$
The domain of this map is minimal non-free, while its range is 
$\aleph_1$-free, hence $\delta$ must be $0$, a contradiction.
$\hfill \square$\\

\medskip

\begin{step-lemma}\label{step}
Let $G = \bigcup\limits_{\alpha \in \omega_1} G_\alpha$
be a filtration of type I and let
$$0 \arr H^0 \arr H_\alpha \overset{h}\arr G_\alpha \arr 0$$
be a short exact sequence with $H^0 \cong G$. If 
$\sigma: G_\alpha \arr H_\alpha$ is a splitting map, then there is an 
extension of this sequence such that $\sigma$ does not extend to a splitting
map $\sigma'$  of the new short exact sequence:

\medskip\hfil\begin{tabular}{rcccccl}
$0 \arr$ & $H^0$         & $\arr$ & $H_\alpha$ & $\overset{h}{\underset{\sigma} 
\rightleftarrows}$ & $G_\alpha$ & $\arr 0$\\
         &$\downarrow$ && $\downarrow$ &                  & $\downarrow$ & \\
$0 \arr$ & $H^0$ & $\arr$ & $H'$ & $\overset{h'}{\underset{\sigma'}
{\overrightarrow{\nleftarrow}}}$ & $G_{\alpha +1}$ & $\arr 0$\\
\end{tabular}\medskip\\
Moreover, the vertical maps in the diagram are inclusions and if 
$G_{\alpha +1}/G_\alpha$ is $n$-free-by-1, then $B'_{\alpha +1}$ is a free 
$R$-module of rank $n$ and 
$$H' = \langle H_\alpha \oplus B'_{\alpha+1},\; y''_{\alpha m} : 
m \in \omega \rangle$$ 
and ${B'_{\alpha+1}}$ is mapped under $h' \mbox{ \it \ mod \ } G_\alpha$ 
onto a free maximal pure $R$-submodule of $G_{\alpha+1} / G_\alpha$.
\end{step-lemma}

\medskip
\noindent
{\bf Proof of the Step-Lemma \ref{step}:} \  We will use special elements
$s_m \in R \quad (m \in \omega  )$
to kill extensions. It will help the reader to pose precise conditions on
the choice of the $s_m$'s only when needed, which will be at the
end of the proof. Readers familiar with such proofs will know that we are 
working to produce a p-adic catastrophe.

First we use the fact that $G' = G_{\alpha+1} / G_\alpha$ is minimal non-free,
say $n$-free-by-$1$. By (\ref{gr}) and (\ref{rr}) we have
$$G' = \langle \bigoplus\limits_{i<n} x_i \; R, \; y_m R, \; m \in \omega 
\rangle$$
with the only relations
$$y_{m+1} p_m = y_m + \sum\limits_{i<n} x_i k_{im} \;\;\; (m \in \omega)$$
and coefficients
$$p_m, k_{im} \in R.$$
\noindent
By the last equations we can find $g_{\alpha m} \in G_\alpha$ and $x_{\alpha i},
y_{\alpha m} \in G_{\alpha + 1}$ such that\\

$(*)$ $\begin{cases}
G_{\alpha+1} = \langle G_\alpha, x_{\alpha i} R, y_{\alpha m} R : i < n, \;\;
m \in \omega \rangle\\
\mbox{with the relations } y_{\alpha m + 1} p_m = y_{\alpha m} + \sum\limits_
{i < n} x_{\alpha i} k_{im} + g_{\alpha m} \; (m \in \omega).
\end{cases}$\\
\bigskip
The action of $\sigma$ is known to us on $G_\alpha$, hence we can 
choose a pure element $ 0 \neq z \in H^0$ and let 
$z_m = z s_m $, 
hence $H^0/zR$ is $\aleph_1$-free by purity of 
$z$ in an $\aleph_1$-free $R$-module. We also choose preimages 
$\overline g_{\alpha m} = g_{\alpha m}\sigma \in H_\alpha$, hence 
$\overline g_{\alpha m}h=g_{\alpha m}$. We are now in the position to apply 
Proposition \ref{embed}. Let 
$$ H'= \langle H_{\alpha}\oplus B'_{\alpha+1}, y''_{\alpha m} :
 m \in \omega \rangle \mbox{ with } B'_{\alpha+1} = 
\bigoplus\limits_{i<n} x''_{\alpha i} R \subseteq H'$$
be the extension given by the proposition with the useful relations
\begin{eqnarray}\label{usef}
y''_{\alpha m+1} p_m = y''_{\alpha m} + \sum\limits_{i<n} 
x''_{\alpha i} k_{im} + z s_m + g_{\alpha m}\sigma
\ \ \ \ (m \in \omega)
\end{eqnarray}
and an extended homomorphism $h': H' \arr G_{\alpha+1} $ with
$$ h'\upharpoonright H_\alpha = h,\ x''_{\alpha i} h'= x_{\alpha i} 
\mbox{ and } y''_{\alpha m} h'= y_{\alpha m} \ (i<n, m \in \omega)$$
such that
$$ \ker h' = H^0 \mbox{ and } \Im h' = G_{\alpha+1}.$$
It remains to show the non-splitting property of the Lemma.

\medskip

Suppose that $\sigma ':G_{\alpha +1} \arr H'$ is an extensions of
$\sigma: G_\alpha \arr H$ such that 
$$\sigma '  h' = id_{G_{\alpha +1}}.$$

Now we want to derive a contradiction when choosing the $s_{m}$'s
accordingly (independent of $\sigma '$ !)
We apply $\sigma '$ to ($*$) and get the equations in $H'$:
\bigskip

\noindent
$$(*\sigma ') \ \ \ \ y_{\alpha m+1} \sigma ' p_m = y_{\alpha m} \sigma ' +
\sum\limits_{i<n} x_{\alpha i} \sigma ' k_{im} + g_{\alpha m} \sigma '.$$
\bigskip
If $d_{\alpha m} = y''_{\alpha m} - y_{\alpha m}\sigma ',$ and 
$ e_{\alpha i} = x''_{\alpha i} - x_{\alpha i} \sigma '$ 
then $ d_{\alpha m} \in H^0$
from 
$$d_{\alpha m}h' = (y''_{\alpha m} - y_{\alpha m}\sigma ')h'
= y''_{\alpha m}h' - y_{\alpha m}\sigma ' h'
= y_{\alpha m} - y_{\alpha m} = 0
\mbox{ and } \ker h' = H^0.$$

\noindent
Similarly we argue with $e_{\alpha i}$ and get
$$ d_{\alpha m}, e_{\alpha i} \in H^0.$$

Subtracting $(*\sigma ')$ from (\ref{usef}) 
leads now to a system of equations in $H^0$.
$$d_{\alpha m+1} p_m = d_{\alpha m} + \sum\limits_{i<n} e_{\alpha i} k_{im} 
+ zs_m.$$
We consider the submodule 
$$W = \langle d_{\alpha m} + zR ,\ e_{\alpha i} + zR \ :i < n, 
m \in \omega \rangle_R \subseteq H^0/zR.$$

\noindent
\noindent
The last displayed equations tell us that W is an epimorphic
image of a minimal non-free $R$-module, hence $0$ or non-free 
of finite rank. On the other hand $H^0/zR$ is $\aleph_1$-free as noted
above, hence $W = 0$  or equivalently
$$\langle d_{\alpha m}, e_{\alpha i} \;:\; m \in \omega, i < n\rangle_R
\subseteq z R \cong R$$
The original equations\\
\begin{eqnarray} \label{cont}
d_{\alpha m+1} p_m = d_{\alpha m} + \sum\limits_{i<n} e_
{\alpha i} k_{im} +s_m 
\end{eqnarray}
still hold, but this time require solutions 
$d_{\alpha m}, e_{\alpha i} \in R.$
We get to an end: just choose rational numbers $s_m \in R$
such that (\ref{cont}) has no solutions. 
The existence of these $s_m$'s follows
from Lemma \ref{nosol}. Finally note that dealing with (\ref{cont}) is 
independent of the particular choices of the extensions of $\sigma$ as 
required in the Lemma.
$\hfill \square$\\

\bigskip
\noindent
{\bf Proof of Proposition 5.1:} Let
$$ H^0 \cong G = \bigcup\limits_{\alpha \in \omega_1} G_\alpha $$
be the module of type I. We must show that $\Ext(G,H^0) \neq 0$ and need a 
non-splitting short exact sequence
\begin{eqnarray}\label{nsp}
0 \arr H^0 \arr H \overset{h}\arr G \arr 0
\end{eqnarray}
which we construct inductively as an ascending, continuous chain of
short exact sequences
$$0 \arr H^0 \arr H_\alpha \overset{h_\alpha}\arr G_\alpha \arr 0$$
with union (\ref{nsp}).
Let
$$0 \arr H^0 \arr H_1 \overset{h_1}\arr G_1 \arr 0$$
be the first step with $G_1$ a free $R$-module of countable rank. By
Observation \ref{obs1} and CH we can enumerate all possible splitting maps
$\sigma :G \arr H$ of extensions $h$ as in (\ref{nsp}) of all $h_1$'s by
$\omega_1$, and let $\{ \sigma_\alpha :G \arr H , \ \alpha \in \omega_1 \}$
be such a list. Using the Step-Lemma \ref{step} and the uniqueness in 
Observation \ref{obs1} we can discard any $\sigma_\alpha$ at stage $\alpha$ 
when constructing
$$0 \arr H^0 \arr H_{\alpha + 1} \overset{h_{\alpha +1}} \arr 
G_{\alpha + 1} \arr 0.$$
The resulting extension (\ref{nsp}) can not split.
$\hfill \square$\\

\section{Splitters Of Type II}

\bigskip

\noindent
An $R$-module $G$ is of type II if $G$ has an $\aleph_1$-filtration
$G = \bigcup\limits_{\alpha \in \omega_1} G_\alpha$ of pure submodules
$G_\alpha$ such that $G/G_\alpha$ is $\aleph_1$-free for all non-limit 
ordinals $\alpha \in \omega_1$, see Section 3. In this section we want to
show our second main

\begin{theorem} \label{thmII}
If $G$ is of type II, then $G$ is a splitter if and only
if $G$ is free over its nucleus $R$.
\end{theorem}

\noindent
{\bf Remark} 
Theorem \ref{thmII} includes that strongly $\aleph_1$-free $R$-modules are 
never splitters, except if trivially the module is free. This was very 
surprising to us.

\bigskip
\noindent
{\bf Proof:} If $S = \{ \alpha \in \omega_1 : G/G_\alpha$ is not
$\aleph_1$-free$\}$, then $S$ is a set of limit ordinals by (\ref{chase}), 
and if $\alpha \in S$ we also may assume that $G_{\alpha+1} /G_\alpha$ is
minimal non-free, compare $\S$3.

\noindent
We get a $\Gamma$- invariant $\Gamma (G)$ defined by $S$ modulo the ideal of
thin sets, see e.g. \cite{EM}. If $\Gamma (G) = 0$, then we find a cub
$C \subseteq \omega_1$, with $ C \cap S = \emptyset$ and $G = \bigcup\limits_
{\alpha \in C} G_\alpha$. Let $\alpha_0 = \min C$. Then $G = G_{\alpha_0} \oplus F$
for some free $R$-module $F$, and $G_{\alpha_0}$ is a countable submodule of
$G$ which must be free over $R$ by Hausen's \cite{Ha} result, see also \cite
{GS}. Hence $G$ is free. Note that the hypothesis of $G$ being $\aleph_1$-free is not
used in this case! If $\Gamma (G) \neq 0$ we want to show that $\Ext (G,G) \neq 0$.
Theorem \ref{thmII} can be rephrased as

\begin{eqnarray} \label{newII}
\mbox{If}\;\; G \mbox{ is of type II, then G is a splitter if and only if } 
\Gamma (G) = 0.
\end{eqnarray}

Now assume that $S$ is stationary in $\omega_1$. 
We want to construct some $H \overset{h} \arr G \arr 0$ with kernel 
$\ker h = H^0,\ G$ isomorphic
to $H^0$ by $\gamma$, which does not split. If 
$G_{\alpha} \gamma = H'_\alpha \ (\alpha \in \omega_1)$, then 
$$H^0 = \bigcup\limits_{\alpha \in \omega_1} H'_\alpha$$

\noindent
is a (canonical) $\aleph_1$-filtration of $H^0$ copied from $G$. First we
pick  elements $z_\alpha \in H^0$ such that $z_\alpha R \cong R$ and
$H^0 / z_\alpha R$ is $\aleph_1$-free, e.g. take any basis element from a 
layer $H'_{\alpha+2} \setminus H'_{\alpha+1}$ of the filtration of $H^0$. Then
we define inductively a continuous chain of short exact sequences
($\alpha \in \omega_1$).

\bigskip

$(\beta)$ $\begin{cases}
0 \arr H^0 \arr H^\beta \overset{h_\beta} \arr G_\beta \arr 0\\
\mbox{countable, free submodules } H_\beta \subseteq _* H^\beta, \mbox{ and
ordinals } \beta < \beta'\leq \omega_1
\end{cases}$
\bigskip

\noindent
subject to various conditions. At the end we want in particular 
$H = \bigcup\limits_{\alpha \in \omega_1} H_\alpha
= \bigcup\limits_{\alpha \in \omega_1} H^\alpha$. 

If $\beta = 0$, then $G_0 = 0$ and we take the zero map 
$h_0 : H^0 \arr G_0 \arr 0$ with kernel $H^0$.

Suppose $(\beta)$ is constructed for all $\beta < \alpha$. If $\alpha$ is a
limit, we take unions $h_\alpha = \bigcup\limits_{\beta < \alpha} h_\beta$,
 $H_\alpha = \bigcup\limits_{\beta < \alpha} H_\beta$
and $H^\alpha = \bigcup\limits_{\beta < \alpha} H^\beta$, assuming that at 
inductive steps sequences extend (naturally) by inclusions. Then visibly
$(\alpha)$ holds. 

We may assume that $(\alpha)$ is known, and we want to
construct $(\alpha + 1)$.\\
If $\alpha \not\in S$, then we extend $(\alpha)$ trivially:\\
Put $H^{\alpha +1} = H^\alpha \oplus F_\alpha$ with $F_\alpha$ a free 
$R$-module of the same rank as the free $R$-module 
$G_{\alpha+1}/G_\alpha$. As $G_{\alpha+1} = F'_\alpha \oplus G_\alpha$, we 
may choose an isomorphism $h': F_\alpha \arr F'_\alpha$ and extend $h_\alpha$ 
to $h_{\alpha+1}$ by $h_{\alpha+1} = h_\alpha \oplus h'$. 
Clearly $\ker h_{\alpha+1} = \ker h_\alpha= H^0$ 
and $\Im h_{\alpha+1} = G_{\alpha+1}$.

If $\alpha \in S$, then we must work. 
We have 
$G_{\alpha+1} / G_\alpha = 
\langle B'_{\alpha+1}, y'_m : m \in \omega \rangle$ from (\ref{gr}) and
(\ref{rr}). Hence
\begin{eqnarray}
G_{\alpha+1} = \langle G_\alpha, B_{\alpha+1}, y_{\alpha m} 
R: m \in \omega \rangle, \
B_{\alpha+1} = \bigoplus\limits_{i<n} x_{\alpha i} R
\end{eqnarray}
with relations 
\begin{eqnarray}  \label{relaII}
y_{\alpha m+1} p_m = y_{\alpha m} + \sum\limits_{i<n}
x_{\alpha i} k_{im} + g_{\alpha m} \ \ \ \ (m \in \omega), 
\end{eqnarray}
where $g_{\alpha m} \in G_\alpha$. 
Let $\bar{B}_{\alpha+1} =
\bigoplus\limits_{i<n} \bar{x}_{\alpha i} R$ be a copy of $B_{\alpha+1}$.
Then we pose the following additional conditions on $(\alpha + 1)$.

\bigskip

\noindent
(a) \quad $H^{\alpha+1} h_{\alpha+1} = H_{\alpha+1} h_{\alpha+1} = G_{\alpha+1}$

\bigskip

\noindent
(b) \quad $\bar{B}_{\alpha+1} \subseteq H_{\alpha+1}$

\bigskip

\noindent
(c) \quad $H^{\alpha+1} / H^\alpha \cong G_{\alpha + 1} / G_\alpha$

\bigskip

\noindent
(d) \quad $H^{\alpha+1} / H_\beta + H'_\gamma \mbox{ is free for all }
\beta \leq \alpha,\ \beta \not\in S \mbox{ and } \gamma \in \omega_1$

\bigskip

\noindent
(e) \quad $H^{\alpha+1} / H_\beta + H'_\gamma \mbox{ is } \aleph_1 \mbox
{-free for all } \beta \leq \alpha, \gamma \in \omega_1$

\bigskip

\noindent
(f) \quad $H'_{\alpha+1} \subseteq H_{\alpha+1} \cap H^0 = H'_{(\alpha+1)'}.$

\bigskip

\noindent
We choose preimages $\bar g_{\alpha m} \in H_\alpha$ such that 
$\bar g_{\alpha m} h_{\alpha} = g_{\alpha m}$ and apply Proposition 
\ref{embed} to define the extension 
$$ H^\alpha \subseteq H^{\alpha+1} = \langle H^\alpha \oplus \bar 
B_{\alpha+1},\ \bar y_{\alpha m}R : m \in \omega \rangle$$
with the relations
\begin{eqnarray} \label{rela}
\bar y_{\alpha m+1} p_m = \bar y_{\alpha m} + \sum\limits_{i<n}
\bar x_{\alpha i} k_{im} + z_{\alpha m} s_{\alpha m} + \bar g_{\alpha m}
\end{eqnarray}
where 
$$\bar B_{\alpha+1} \cong \bigoplus\limits_{i<n} x'_i \;R$$ 
as required in (b). 
Similarly, by Proposition \ref{embed} the map $h_\alpha$ extends to an 
epimorphism
$h_{\alpha+1} : H^{\alpha+1} \arr G_{\alpha+1}$. 
It is now easy to check that (c) holds and it is also easy to see that 
$\ker h_{\alpha+1} = H^0$.
Next we extend $H_\alpha \subset H_{\alpha+1}$ carefully such that 
$(\alpha + 1)$, (a), (b), (d), (e) and (f) hold.

\medskip
\noindent
$\Im h_{\alpha+1} = G_{\alpha+1}$ is a countable module of the $\aleph_1$-
free $R$-module $G$, hence free and $h_{\alpha+1}$ must split. There is a
splitting map
$$\sigma : G_{\alpha+1} \arr H^{\alpha+1}\;\;\mbox{such that}\; \sigma h_{\alpha+1}
= id_{G_{\alpha+1}},$$
hence $H^{\alpha+1} = H^0 \oplus G_{\alpha+1} \sigma$.  Let
$$\pi_0 : H^{\alpha+1} \arr H^0,\ \pi_1 : H^{\alpha+1} \arr
G_{\alpha+1} \sigma$$
be the canonical projections with 
$\pi_0 + \pi_1 = id_{H^{\alpha+1}}$. Recall
that $\bar B_{\alpha+1} \subseteq H^{\alpha+1}$. Choose $\beta \in \omega_1$
large enough such that $\beta \not\in S, \; \beta \geq \alpha', \alpha+1$ and
$(\bar{B}_{\alpha+1} + H_\alpha) \pi_0 \subseteq H'_\beta$. This is easy
because $\omega_1 \setminus S$ is unbounded and $(\bar{B}_{\alpha+1} + H_\alpha)
\pi_0$ is countable. Put $H_{\alpha+1} = H'_\beta \oplus (G_{\alpha+1} \sigma)$
and $\beta = (\alpha+1)'$. Note that 
$$G_{\alpha+1} = H^{\alpha+1} h_{\alpha+1}
\supseteq H_{\alpha+1} h_{\alpha+1} \supseteq G_{\alpha+1} \sigma h_{\alpha+1}
= G_{\alpha+1}$$
and (a) follows.\\

\noindent
If $\gamma \leq (\alpha+1)'$, then $H_{\alpha+1} + H_\gamma = H_{\alpha+1}$
and if $\gamma \geq (\alpha+1)'$, then $H_{\alpha+1} + H'_\gamma=H'_\gamma \oplus
G_{\alpha+1} \sigma$ and (d) follows. We see immediately $H_\alpha \pi_1
\subseteq G_{\alpha+1} \sigma$ and 
$H_\alpha \pi_0 \subseteq H'_{(\alpha+1)'}$, 
hence
$$H_\alpha \subseteq H_\alpha \pi_0 + H_\alpha \pi_1 \subseteq H_{\alpha+1}$$

\noindent
and $(\alpha+1)$ holds; similarly $\bar{B}_{\alpha+1} \subseteq H_{\alpha+1}$
for (c). From $H'_{(\alpha+1)'} \subseteq H^0$ and the modular law we have
$H_{\alpha+1} \cap H^0 = H'_{(\alpha+1)'} \oplus (G_{\alpha+1} \sigma \cap
H^0) = H'_{(\alpha+1)'}$ and (f) holds.

\bigskip

\noindent
Finally we choose $H = \bigcup\limits_{\alpha \in \omega_1} H^\alpha ,
h = \bigcup\limits_{\alpha \in \omega_1} h_\alpha$ and
\begin{eqnarray} \label{nonsp}
 0 \arr H^0 \arr H \overset{h} \arr G \arr 0
\end{eqnarray}
is established and it remains to show that (\ref{nonsp}) does not split. 
Suppose for contradiction that $\sigma : G \arr H$ is a splitting map 
for $h$. We have
$H = \bigcup\limits_{\alpha \in \omega_1} H_\alpha$ and 
$G = \bigcup\limits_{\alpha \in \omega_1} G_\alpha$ $\aleph_1$-filtration. 
Using the above properties of the $H_\alpha$'s, it follows by a back and 
forth argument that 
$$E = \{\alpha \in \omega_1: H_\alpha \cap H^0 = H'_\alpha, \; 
G_\alpha \sigma \subseteq H_\alpha\}$$
\noindent
is a cub. On the other hand $S$ is stationary in $\omega_1$ and we find
$$\alpha \in S \cap E.$$

\noindent
From (\ref{relaII}) and (\ref{rela}) we have
$$y_{\alpha m+1} \sigma p_m = y_{\alpha m} \sigma + \sum\limits_{i<n} 
x_{\alpha i} \sigma k_{im} + g_{\alpha m} \sigma$$

\noindent
and
$$\bar{y}_{\alpha m+1} p_m = \bar{y}_{\alpha m} + \sum\limits_{i<n} \bar{x}_
{\alpha i} k_{im} + z_\alpha s_{\alpha m} + \bar{g}_{\alpha m}\;\;\mbox{with}\;\;
g_{\alpha m} \sigma, \; \bar{g}_{\alpha m} \in H_\alpha\;\;\mbox{and}\;\;
G_\alpha \sigma \subseteq H_\alpha.$$
\bigskip
\noindent
Put $d_{\alpha m} = \bar{y}_{\alpha m} - y_{\alpha m} \sigma,\ 
f_{\alpha m} = \bar{g}_{\alpha m} - g_{\alpha m} \sigma, \
e_{\alpha i} = \bar{x}_{\alpha i} - x_{\alpha i} \sigma$ and notice that
$d_{\alpha m} h = e_{\alpha i} h = f_{\alpha m}h = 0$, 
hence $d_{\alpha m}, e_{\alpha i}, f_{\alpha m} \in H^0.$\\

\bigskip
\noindent
Subtracting the last displayed equations we get\\

\noindent
(j) \quad
$d_{\alpha m+1} p_m = d_{\alpha m} + \sum\limits_{i<n} e_{\alpha i}
k_{im} + f_{\alpha m} + z_\alpha s_{\alpha m} \mbox{ in } H^0.$\\

Recall that 
$f_{\alpha m} \in H_\alpha \cap H^0 \subseteq H'_{\alpha'}$ by
(f) and  
{ \it mod \ } $T = H'_{\alpha'} + z_\alpha R$ 
the equations (j) say that $W =
\langle d_{\alpha m}, e_{\alpha i} : i < n, m \in \omega
\rangle + T/T$ is either minimal non-free or 0. On the other hand
$H^0/T$ is $\aleph_1$-free, hence $W = 0$ and $\langle d_{\alpha m},
e_{\alpha i}, z_\alpha \rangle \subseteq H'_\alpha + z_\alpha R$.
Recall from (e) that $H^0 / H'_\alpha$ is ${\aleph_1}$-free. Hence (j) turns
into
$$d_{\alpha m+1} p_m \equiv d_{\alpha m} + \sum\limits_
{i<n} e_{\alpha i} k_{im} + z_\alpha s_{\alpha m} \mod H'_\alpha.$$

\noindent
Using $\aleph_1$-freeness of $H^0/z_\alpha R$ these equations tell us that
we must have solutions $d_{\alpha m}, e_{\alpha i} \in R$ for\\

\bigskip
\noindent
(k) \quad $d_{\alpha m+1} p_m = d_{\alpha m} + \sum\limits_{i<n} e_{\alpha i}
k_{im} + s_{\alpha m}.$\\

\medskip
\noindent
In Lemma \ref{nosol} we selected particular $s_{\alpha m}$'s in $R$ 
such that (k) has no solution in $R$. Now we are ready to make this choice 
which we should have done right at the beginning of the proof and 
hence derive a contradiction; we conclude $\Ext (G,G) \neq 0$.
$\hfill \square$\\

\bigskip

\noindent
From Theorem \ref{thmII} we see that non-free but strongly 
$\aleph_1$-free abelian groups are never splitters. 
We find this very surprising. Particular groups like the Griffith-group 
$\GG$ below which is a Whitehead group $(\Ext (\GG, \Z) = 0)$ under 
Martin's axiom and $\neg CH$ is not a splitter. Recall a nice and easy 
construction of  $\GG$ which is sometimes Whitehead but always fails to be a
splitter in general.

\bigskip
\noindent
Let $P = \Z^{\aleph_1} = \prod\limits_{\alpha \in \aleph_1} \alpha \Z$ the
cartesian product of $\Z$. If $\lambda \in \aleph_1$ is a limit ordinal 
choose an order preserving map $\delta_\lambda : \omega \ra \lambda$ 
with sup$(\omega \delta_\lambda) = \lambda$. Then, along this ladder 
system we define branch-elements
$$c_{\lambda n} = \sum\limits_{i \geq n} (i \delta_\lambda) \frac{i!}{n!}$$
which are a `divisibility chain' of $c_{\lambda_{0}}$ modulo 
$\bigoplus\limits_{\alpha \in \aleph_1} \alpha \Z$, hence 
$$\GG = \langle \bigoplus\limits_{\alpha \in \aleph_1} \alpha \Z, 
c_{\lambda n} : \lambda \in \aleph_1, \  \lambda 
\mbox{ a limit ordinal, } n \in \omega\rangle$$ 
is a pure subgroup
of $P$. We see that $|\GG| = \aleph_1$ and $\GG$ is $\aleph_1$- free by
$\aleph_1$-freeness of $P$; see \cite{Fu} (Vol 1, p 94, Theorem 19.2).
Moreover $\Gamma_{G} \neq 0$ because $G_\beta = \GG \cap \prod\limits_{\alpha < \beta}
\alpha \Z$ $(\beta\in\omega_1)$ is an $\aleph_1$-filtration of $\GG$ with $G_{\lambda+1} / G_\lambda$
divisible for all limit ordinals $\lambda$. Hence $\GG$ is not free. It is
easy to check that $\GG$ is $\aleph_1$-separable, hence strongly $\aleph_1$-free;
see also \cite{EM}, p. 183, Theorem 1.3.
\bigskip

\section{Splitters Of Type III}

If $G$ is of type III then we recall from Section 3 that $G = \bigcup\limits_
{\alpha \in \omega_1} G_\alpha ,\ G_0 = 0$ with (3.4) -- (3.6) and $G_1 =
\bigcup\limits_{j \in \omega_1} G_{0j}$ and $\{G_{0j} : j \in \omega_1\}$
is an $\aleph_1$-filtration of pure submodules $G_{0j}$ such that each
$G_{0 j+1} / G_{0j}$ is minimal non-free. Here we will show:

\begin{theorem}  \label{thmIII}
Modules of type III are not splitters.
\end{theorem}

\noindent
{\bf Proof:} Let $G' = \bigcup\limits_{\alpha \in \omega_1} G'_\alpha$
be an isomorphic copy of $G$ taking $G_\alpha$ to $G'_\alpha$, and choose a
sequence of elements $z_\alpha \in G'_{\alpha + 2} \ (\alpha \in \omega_1)$
such that 
$$G'_{\alpha+1} \cap z_\alpha R = 0  \mbox{ and } 
G'_{\alpha + 2}/G'_{\alpha+1}\oplus z_\alpha R \mbox{ is } 
\aleph_1\mbox{-free. }$$ 
This is possible by (III).     

\bigskip

\noindent
By a basic observation from Section 1 it is enough to show that $\Ext (G_1, G)
\neq 0$. Inductively we will construct a non-trivial element in $\Ext (G_1, G)$.
We consider the following diagram\\

\medskip\hfil\begin{tabular}{rcccccl}

(0)\quad\quad $0 \arr$ & $G'$ & $\arr$ & $H^0$ & $\overset{h_0}\arr$ & $G_{00} = 0$ & $\arr$ $0$\\

         &$\|$ && $\downarrow$ && $\downarrow$ & \\

$(\beta)$\quad\quad $0 \arr$ & $G'$ & $\arr$ & $H^\beta$ & $\overset{h_\beta}\arr$ & $G_{0 \beta}$ & $\arr$  $0$\\

         &$\|$ && $\downarrow$ && $\downarrow$ & \\

$(\omega_1)$\quad\quad  $0 \arr$ & $G'$ & $\arr$ & $H$ & $\overset{h} \arr$ & $G_1$ & $\arr$  $0$

\end{tabular}\medskip

\bigskip

\noindent
The first row is the trivial extension with $G' = H^0$ and $h_0 = 0$. Vertical
maps and maps between $G'$ and $H'$s are inclusions. The sequences $(\beta)$
are increasing continuous and suppose $(\beta)$ is constructed for all
$\beta < \alpha$. Then $h_\alpha = \bigcup\limits_{\beta < \alpha} h_\beta$
and \\
$$0 \arr G' \arr \bigcup\limits_{\beta < \alpha} H^\beta \overset{h_\alpha} \arr
G_{0 \alpha} \arr 0$$

\noindent
if $\alpha$ is a limit. Next we want to construct $(\alpha + 1)$ from
$(\alpha)$ and recall that $G' = G_{0 \alpha + 1} / G_{0 \alpha}$ is minimal
non-free generated as in (\ref{gr}), (\ref{rr}). We can write

\begin{eqnarray}
G_{0 \alpha + 1} = \langle G_{0 \alpha}, \ B_{\alpha + 1}, y_{\alpha m} R 
: m \in\omega \rangle, \ \ B_{\alpha+1} = 
\bigoplus\limits_{i<n} x_{\alpha i}R
\end{eqnarray}

with relations 

\begin{eqnarray} \label{relaIII}
y_{\alpha m+1} p_m = y_{\alpha m} + \sum\limits_{i<n} x_{\alpha i} k_{im} 
+ g_{\alpha m} \; (m \in \omega), \; g_{\alpha m} \in  G_{0 \alpha}. 
\end{eqnarray}

\noindent
Then we define 
$$h_{\alpha+1} : H^{\alpha+1} \arr G_{0 \alpha + 1}\arr 0.$$

\noindent
by Proposition \ref{embed}. Hence
$$H^{\alpha+1} = \langle H^\alpha \oplus \bar{B}_{\alpha+1} ,\ 
\bar{y}_{\alpha m}R : m \in \omega \rangle$$
\noindent
has the relations

\begin{eqnarray}\label{relH}
\bar{y}_{\alpha m+1} p_m = \bar{y}_{\alpha m} + \sum\limits_{i<n}
\bar{x}_{\alpha i} k_{im} + z_\alpha s_{\alpha m} + \bar{g}_{\alpha n} \;\;
(m \in \omega),
\end{eqnarray}

where the $s_{\alpha m} \in R$ will be specified later on, and 
$\bar{g}_{\alpha m} \in H^\alpha$. 

Suppose that $(\omega_1)$ splits and consequently $\sigma : G_1 \arr H$
is a splitting map for $h$. Then let 
$$d_{\beta m} = y_{\beta m} - \bar{y}_{\beta m} ,
\ e_{\beta i} = x_{\beta i} \sigma - \bar{x}_{\beta i} \mbox{ and }
f_{\beta m} = g_{\beta m} \sigma - \bar{g}_{\beta m}.$$ 
From
splitting we get again 
$$d_{\beta m},\ e_{\beta i}, \ f_{\beta m} \in G' \ \ \  (\beta \in 
\omega_1, m \in \omega, i<n).$$ 
Using
$G' = \bigcup\limits_{\alpha \in \omega_1} 
G'_\alpha$ for $\beta \in \omega_1$
we find $\alpha \in \omega_1$ such that 

$$d_{\beta m}, \ e_{\beta i}, \ f_{\beta m} \in G'_\alpha \mbox{ for all } 
i<n, m \in \omega.$$
Consider a map $\tau : \omega_1 \arr \omega_1$ taking any 
$\alpha \in \omega_1$ to
$$\tau (\alpha) = \mbox{ min } 
\{\beta \in \omega_1; \beta  \mbox{ a limit ordinal, }
d_{\alpha m},\ e_{\alpha i},\ f_{\alpha m} \in G'_\beta, \alpha \leq \beta, 
m \in \omega , i < n \}$$
and note that 
$C = \{ \alpha \in \omega_1 : \tau (\alpha)= \alpha \}$ 
is a cub in
$\omega_1$ and a subset of
$$E = \{ \alpha \in \omega: \alpha  \mbox{ a limit ordinal, }  d_{\beta m}, 
\ e_{\beta i},\ f_{\beta m} \in G'_\alpha \mbox{ for all } \beta < \alpha, 
i < n, m \in \omega\}.$$
Hence $E$ is a cub in $\omega_1$. Next we apply $\sigma$ to (\ref{relaIII}) 
and subtract (\ref{relH}). Hence we get a system of equations in 
$H^{\alpha+1}$.
\begin{eqnarray}\label{diff}
d_{\alpha m+1} p_m = d_{\alpha m} + \sum\limits_
{i<n} e_{\alpha i} k_{im} + f_{\alpha m} + z_\alpha
s_{\alpha m} \ \ (m\in\omega)
\end{eqnarray}
If $\alpha \in E$, then {\it modulo} $H^\alpha$ the equations 
(\ref{diff}) turn into
$$d_{\alpha m+1} p_m \equiv d_{\alpha m} + \sum\limits_
{i<n} e_{\alpha i} k_{im} + z_\alpha s_{\alpha m} \;\; (m \in \omega)$$
and modulo $z_\alpha R$ an earlier argument and $\aleph_1$-freeness of
$G' / G_\alpha \oplus z_\alpha R$ show that the last equation requires
solutions $d_{\alpha m}, e_{\alpha i} \in R$ for
$$d_{\alpha m+1} p_m = d_{\alpha m} + \sum\limits_{i<n} e_{\alpha i} k_{im}
+ s_{\alpha m} \;\; (m \in \omega)$$
By a special choice of $s_{\alpha m}$'s in Lemma \ref{nosol} this is now 
excluded, a contradiction. Hence $(\omega_1)$ has no splitting map and 
Theorem \ref{thmIII} 
follows.$\hfill \square$\\

\section{Appendix: Splitters Of Type I Under $2^{\aleph_0} < 2^{\aleph_1}$}

\noindent
In Section 5 we have seen a proof that CH implies  modules 
of type I are never splitters. 
A slight variation but some what technical modification of the proof, 
shows that this result can be extended 
to WCH that is $2^{\aleph_0} < 2^{\aleph_1}$. Due to Section 5 
this is not needed for the main result of this paper dealing with 
modules of cardinality  $\aleph_1$ but it will be interesting when 
passing to cardinals  $> \aleph_1$.
We outline the main steps, their proofs are suggested by the 
 proofs in Section 5. 
\begin{theorem} \;\; {\rm ( ZFC $ + \ 2^{\aleph_0} < 2^{\aleph_1}$)}
Modules of type I are not splitters. 
\end{theorem}
\begin{step-lemma} \label{wstep}
Let $G = \bigcup\limits_{\alpha \in \omega_1} G_\alpha$
be a filtration of type I and let
$$0 \arr K \arr H_\alpha \overset{h}\arr G_\alpha \arr 0$$
be a short exact sequence with some $z \in K$ such that $z R \cong R$ and $K,
K / z R$ are $\aleph_1$-free. Then there are two commuting diagrams $(\epsilon = 0,1)$

\medskip\hfil\begin{tabular}{rcccccl}
$0 \arr$ & $K$         & $\arr$ & $H_\alpha$ & $\overset{h} \arr$ & $G_\alpha$ & $\arr 0$\\
         &$\downarrow$ && $\downarrow$ &                  & $\downarrow$ & \\
$0 \arr$ & $K_\epsilon$ & $\arr$ & $H_\epsilon$ & $\overset{{h_\epsilon}} \arr$ & $G_{\alpha+1}$ & $\arr 0$\\
\end{tabular}\medskip

\noindent
with vertical maps inclusions such that any third row with $H'_\epsilon$
$\aleph_1$-free,

$$0 \arr K'_\epsilon \arr H'_\epsilon \overset{{h'_\epsilon}} \arr G_\beta \arr 0$$
and any splitting map $\sigma$ of $h$ cannot have two splitting extensions
$\sigma_\epsilon$ of $h'_\epsilon:$

\medskip\hfil\begin{tabular}{rcccccl}
$0 \arr$ & $K$         & $\arr$ & $H_\alpha$ & $\overset{h}{\underset{\sigma} \rightleftarrows}$
 & $G_\alpha$ & $\arr 0$\\
         &$\downarrow$ && $\downarrow$ &                  & $\downarrow$ & \\
$0 \arr$ & $K_\epsilon$ & $\arr$ & $H_\epsilon$ & $\overset{{h_\epsilon}} \arr$ & $G_{\alpha+1}$ & $\arr 0$\\
         &$\downarrow$ && $\downarrow$ &                  & $\downarrow$ & \\
$0 \arr$ & $K'_\epsilon$ & $\arr$ & $H'_\epsilon$
& $\overset{h_\epsilon'}{\underset{\sigma_\epsilon}{\overrightarrow{\nleftarrow}}}$ & $G_\beta$ & $\arr 0$\\
\end{tabular}

\medskip
\noindent
Moreover $H_\epsilon =
\langle H_\alpha \bigoplus \overline B^\epsilon_{\alpha+1},\;
y^\epsilon_{\alpha m} : m \in \omega \rangle$ and
$\overline B^\epsilon_{\alpha+1}$ is
mapped under $h_\epsilon \mod G_\alpha$ onto a free maximal $R$-submodule of
$G_{\alpha+1} / G_\alpha$ cf. (\ref{gr}).
\end{step-lemma}

\medskip

\noindent
\begin{definition}  
If such an extension $\sigma_\epsilon$ as in (\ref{wstep}) exists
for some $\epsilon \in \{0,1\}$ we say that $\sigma$ {\bf splits over} 
$(H_\epsilon, h_\epsilon).$
\end{definition}

\bigskip

\noindent
{\bf Proof of Lemma \ref{wstep}:}  Compare the proof of
the Step-Lemma \ref{step} but note that at the end you must take once more
differences of the elements $d_{\alpha m}, e_{\alpha i}$ for $\epsilon =1$
and $\epsilon = 0$ respectively. Then we are able to apply Lemma 
\ref{nosol} to get a contradiction from splitting. $\hfill \square$\\

\noindent
We then apply the Step-Lemma and weak diamond $\Phi_{\omega_1}$ to construct 
a short exact sequence
$$0 \arr H^0 \arr H \overset{h}\arr G \arr 0.$$

Let $\gamma: G\arr H^0$ be a fixed isomorphism. Later we will use 
consequences of $\Phi_{\omega_1}$
to show that $h$ does not split.\\

\noindent
{\bf Proof of Theorem 8.1:}\;\; If $H'_\alpha = G_\alpha \gamma$, then $H^0 = \bigcup
\limits_{\alpha \in \omega_{1}} H'_\alpha$ is an $\aleph_1$-filtration if
$G = \bigcup\limits_{\alpha \in \omega_{1}} G_\alpha$ is the given filtration
of type I.\\

\medskip
\noindent
Let $T = {}^{\omega_{1}>} 2$ be the tree of all branches $\eta : \alpha \arr 2$
for some $\alpha \in \omega_1$. We call $\alpha = l (\eta)$ the length of
$\eta.$ Branches are ordered as usually, hence $\eta < \eta'$ if $\eta'
\upharpoonright \Dom \; \eta = \eta.$ The empty set $\emptyset$ is the bottom
element of the tree. If $\eta \in T$, then we construct triples
$$(H_\eta, H^\eta, h_\eta)$$
of $R$-modules $H_\eta \subseteq H^\eta$ with $H_\eta$ free of countable rank
and a homomorphism
$$h_\eta : H^\eta \arr G$$
subject to various natural conditions.

\bigskip
\noindent
(i) \quad $\begin{cases}
H^\emptyset = H^0, H_\emptyset = 0 \;\;\mbox{and}\;\;
h_\emptyset = 0, \;\;\mbox{hence}\\
0 \arr H^0 \arr H_\emptyset \overset{{h_\emptyset}} \arr G_0 = 0 \arr 0
\text{ is short exact.}
\end{cases}$\\

\bigskip

\noindent
(ii) \quad $\begin{cases}
\mbox{If}\;\; \eta < \eta', \;\;\mbox{then}\;\;
(H_\eta, H^\eta, h_\eta) \subseteq (H_{\eta'}, H^{\eta'}, h_{\eta'}),\mbox{ i.e.}\\
H_\eta \subseteq H_{\eta'}, H^\eta \subseteq H^{\eta'} \; \mbox{and}\; h_\eta
\subseteq h_{\eta'}.
\end{cases}$\\

\bigskip

\noindent
(iii) \quad $H^0 = \ker h_\eta \; \mbox{and}\; \Im h_\eta = G_{l(\eta)} = H_\eta h_\eta.$

\bigskip

\noindent
(iv) \quad $\begin{cases}
\text{If }\alpha \in \omega_1 \text{ is a limit and }\eta \in {}^{\alpha}2,
\text{ then we take unions}\\
(H_\eta, H^\eta, h_\eta) = \bigcup\limits_{\beta < \alpha} \; (H_{\eta 
\upharpoonright \beta}, H^{\eta \upharpoonright \beta}, h_{\eta \upharpoonright \beta})
= (\bigcup\limits_{\beta < \alpha} H_{\eta \upharpoonright \beta}, \bigcup
\limits_{\beta < \alpha} H^{\eta \upharpoonright \beta}, \bigcup\limits_
{\beta < \alpha} h_{\eta \upharpoonright \beta}).
\end{cases}$\\

\bigskip

\noindent
If $l (\eta) = \alpha $  we put further
restrictions on those triples. In this case $G' = G_{\alpha + 1} / G_\alpha$ is
minimal non-free, and $G'$ can be represented by (\ref{gr}), (\ref{rr}). There are
elements $g_{\alpha m} \in G_\alpha, x_{\alpha i}, y_{\alpha m} \in 
G_{\alpha + 1}$ with 
$$G_{\alpha+1} = \langle G_\alpha, B_{\alpha+1}, y_{\alpha m} R \;: \; m \in
\omega \rangle, \;\; B_{\alpha+1} = \bigoplus\limits_{i<n} x_{\alpha i} \; R$$
and relations
$$y_{\alpha m+1} p_m = y_{\alpha m} + \sum\limits_{i<n} x_{\alpha i} k_{im} +
g_{\alpha m}.$$
We choose an isomorphic copy $\overline B_{\alpha+1} = 
\bigoplus\limits_{i<n}
\overline x_\alpha R$ of $B_{\alpha+1}$ and now continue defining the tree with
triples.\\

\medskip
\noindent
If $\epsilon \in \{0,1\}$, then we require more from $(H_{\eta\hat{~}<\epsilon>},
H^{\eta\hat{~}<\epsilon>}, h_{\eta\hat{~}<\epsilon>})$.

\begin{align*}
\text{(S\;i)} \quad \quad & H_{\eta} \oplus \overline B_{\alpha+1} \subseteq H_{\eta \hat{~} <\epsilon>}\\
\text{(S\;ii)} \quad\quad & H_{\eta'}  + H'_\beta \subseteq_* H^{\eta \hat{~} <\epsilon>} \mbox{ for all }
\eta' \in {}^{\alpha \geq}2, \; \mbox{and}\; \beta \in \omega_1.\\
\text{(S\;iii)} \quad\quad & H'_\alpha \subseteq H_\eta \cap H^0 \subseteq H_{\alpha '} \mbox
{ for some } \alpha' \in [\alpha, \omega_1)
\end{align*}

\noindent
Note that $H_\eta \cap H^0 = \ker (h_\eta \upharpoonright H_\eta).$\\

\medskip

\noindent
(S\;iv) The crucial condition:\\

\noindent
Suppose $\sigma:  G_\alpha \arr H_\eta$ is a homomorphism
extending to $\sigma_\epsilon: G_{\alpha+1} \arr H^{\eta \hat{~} <\epsilon>}$,
then not both of them can be splitting maps over $(H_{\eta \hat{~} <\epsilon>},
h_{\eta \hat{~} <\epsilon>})$ for $\epsilon = 0,1.$

\bigskip

\noindent
Before we begin with the inductive construction, we observe from (iii) that
$G_{\alpha+1} / G_\alpha \cong H_{\eta \hat{~} <\epsilon>} / H_\eta \cong 
H^{\eta \hat{~} <\epsilon>} / H^\eta$ for $\Dom \eta = \alpha.$\\

\medskip

\noindent
If $\eta \in {}^{\omega_1}2$ and $H = H (\eta) = \bigcup\limits_{\alpha < \omega_1}
H^{\eta \upharpoonright \alpha},$ then $H'= \bigcup\limits_{\alpha < \omega_1}
H_{\eta \upharpoonright \alpha} \subseteq H$ from (ii). (S\;iii) ensures
$H^0 \subseteq H'$ and from (iii) we get $H/H_0 = H'/H_0$, hence $H'=H.$ This
will show that 
$$(\eta 1) \quad\quad H (\eta) = \bigcup\limits_{\alpha < \omega_1}
H_{\eta \upharpoonright \alpha} = \bigcup\limits_{\alpha < \omega_1} H^{\eta
\upharpoonright \alpha} \quad\quad (\eta \in {}^{\omega_1}2)$$

\medskip

\noindent
Similarly 
$h (\eta) = \bigcup\limits_{\alpha < \omega_1} 
h_{\eta \upharpoonright \alpha}$
is a well-defined homomorphism $h (\eta) : H(\eta) \arr G$ by (ii), 
it is onto 
with kernel $H^0$ by (iii), hence 
$$(\eta 2) \quad\quad 0 \arr H^0 \arr H (\eta) \overset{h(\eta)} 
\arr G \arr 0 \quad\quad  (\eta \in {}^{\omega_1}2).$$
\noindent
Condition $(\eta 1)$ provides an $\aleph_1$-
filtration used to apply weak diamond $\Phi$ for showing that 
$(\eta 2)$ does not split for some $\eta.$

\medskip

\noindent
Next we will show that the tree with triples exists. This will follow by
induction along the length $\alpha$ branches $\eta \in  {}^\alpha 2$. The case
$\alpha = 0$ is (i) and already established. Suppose the construction is
completed for all $\beta < \alpha$ and $\alpha < \omega_1$ is a limit
ordinal. For $\eta \in {}^\alpha 2$ we define $(H_\eta, H^\eta, h_\eta)$ as in
(iv) and it is easy to verify that all conditions hold, notably (S\;iii),
because we take only countable unions. We come to the inductive step 
constructing $(H_{\eta \hat{~} <\epsilon>}, H^{\eta \hat{~} <\epsilon>}, h_{\eta \hat{~} <\epsilon>})$
from $(H_\eta, H^\eta, h_\eta)$ for $\alpha = \Dom \eta$.

\bigskip

\noindent
First we adopt the Step-Lemma for $K = H^\eta \cap H^0 = \ker h_\eta, H_\alpha =
H^\eta, H_\epsilon = H^{\eta \hat{~} <\epsilon>}, h = h_\eta, h_\epsilon = 
h_{\eta \hat{~} <\epsilon>}$ and note that the needed element $z$ exists because
$H^\eta \cap H^0 \subseteq H'_\alpha$ is free. We must still define 
$H_{\eta \hat{~} <\epsilon>} \supseteq H_\eta$ carefully satisfying (ii), (S\;i) -
(S\;iii) and the last equality in (iii): Write again $h_\epsilon$ for
$h_{\eta \hat{~} <\epsilon>}$ and $H^\epsilon$ for $H^{\eta \hat{~} <\epsilon>}.$
We know that $\Im h_\epsilon = G_{\alpha+1}$ is a countable submodule of the
$\aleph_1$-free module $G$, hence free and $h_\epsilon$ must split. There is
a splitting map $\varphi_\epsilon : G_{\alpha+1} \arr H^\epsilon$ such that
$\varphi_\epsilon h_\epsilon = id_{G \alpha+1}$, hence
$$H^\epsilon = H^0 \oplus (G_{\alpha+1} \varphi_\epsilon)$$
\noindent
from the first part of (iii). Let $\pi^\epsilon_0 : H^\epsilon \arr H^0$
and $\pi^\epsilon_1 : H^\epsilon \arr H^0$ be the canonical projections, hence
$\pi^\epsilon_0 + \pi^\epsilon_1 = id_{H^\epsilon}$.

\bigskip

\noindent
Choose $\beta = (\alpha+1)' < \omega_1$ large enough such that
$\overline B_{\alpha+1} \pi_0 \cup H_\eta \pi_0 \subseteq H'_\beta,$
where $\overline B_{\alpha+1}$ is taken from the Step-Lemma. We can choose
$\beta$ because $\overline B_{\alpha+1}$ and $H_\eta$ are countable. Put
$$H_{\eta \hat{~} <\epsilon>} = H'_\beta \oplus (G_{\alpha+1} \varphi_\epsilon),$$

\noindent
hence by the known half of (iii)
$$G_{\alpha+1} = H^\epsilon h_\epsilon \supseteq H_{\eta \hat{~} <\epsilon>}
h_\epsilon \supseteq G_{\alpha+1} \varphi_\epsilon h_\epsilon = G_{\alpha+1}$$

\noindent
and the other half of (iii) follows.\\
\noindent
If $\gamma \leq \beta$, then $H_{\eta \hat{~} <\epsilon>} + H'_\gamma =
H_{\eta \hat{~} <\epsilon>}$ and if $\gamma \geq \beta$, then 
$H_{\eta \hat{~} <\epsilon>} + H'_\gamma =
H'_\gamma \oplus (G_{\alpha+1} \varphi_\epsilon)$ with quotient 
$H^{\eta \hat{~} <\epsilon>} / H'_\gamma \oplus (G_{\alpha+1} \varphi_\epsilon)
\cong H^0 / H'_\gamma$ which shows (S\;ii). Trivially $H_\eta \pi_1 \subseteq
(G_{\alpha+1} \varphi_\epsilon)$ and $H_\eta \pi_0 \subseteq H'_\beta$ by the
choice of $\beta$, hence
$$H_\eta = H_\eta id_{H^\epsilon} \subseteq H_\eta \pi_0 + H_\eta \pi_1
\subseteq H_{\eta \hat{~} <\epsilon>}$$
and (ii) holds. Similarly $\overline B_{\alpha+1} \subseteq H_{\eta \hat{~} < \epsilon >}$
and (S\;i) is shown. From $H'_\beta \subseteq H^0$ and the modular law we
have
$$H_{\eta \hat{~} <\epsilon>} \cap H^0 = H'_\beta \oplus (G_{\alpha+1} \varphi_
\epsilon \cap H^0) = H'_\beta$$
and (S\;iii) holds.\\

\bigskip

The construction of the tree with triples is complete. We are ready to use the
weak diamond $\Phi_{\aleph_1}(S)$ to show that $G$ is not a splitter.

\bigskip

\noindent
We will use $\Phi_{\omega_1}(S)$ as stated in Eklof, Mekler 
\cite[p. 143,  
Lemma 1.7]{EM} and note that $G = \bigcup\limits_{\alpha < \omega_1} G_\alpha ,
H (\eta) = \bigcup\limits_{\alpha < \omega_1} H_{\eta \upharpoonright \alpha}$
are $\aleph_1$-filtrations. $\Phi_{\omega_1}(S)$ must tell us which 
$\eta \in {}^{\omega_1} 2$ we should pick. We define a partition $P$ so that
for $\alpha \in S$ a homomorphism $\sigma : G_\alpha \arr H_\eta$ $(\eta \in
{}^\alpha 2)$ has value $P_\alpha (\sigma) = 0$ if and only if $\sigma$ does
{\bf not} split over $(H_{\eta \hat{~} 0}, h_{\eta \hat {~} 0})$. By 
the Step-Lemma build into the construction, we observe that

\bigskip

$(a)$ $
\mbox{if}\; P (\sigma) = 1, \;\mbox{then}\; \sigma \mbox{ cannot split over }
(H_{\eta \hat{~} 1}, h_{\eta \hat {~}1})
$
 
\bigskip
\noindent
The prediction principle finds us a branch $\eta \in {}^{\omega_1} 2$ with the
$\Phi$-property

\bigskip

$(b)$ $\begin{cases}
\mbox{If}\; \sigma : G \arr H \mbox{ is any map, then}\\
S'=\{\alpha \in S : P_\alpha (\sigma \upharpoonright G_\alpha) = \eta (\alpha)\}
\mbox{ is stationary in } \omega_1.
\end{cases}$

\bigskip

\noindent
We pick that branch and build $H = H (\eta)$ and $h = h(\eta)$ accordingly,
hence\\ 
$0 \arr H^0 \arr H \overset{h} \arr G \arr 0$ is short exact. After the branch
is fixed we let $H^{\eta \upharpoonright \alpha} = H^\alpha, h \upharpoonright
H^\alpha = h_\alpha$ and $H_{\eta \upharpoonright \alpha} = H_\alpha$. Now we
claim that the last sequence does not split. Suppose to the contrary that 
$\sigma : G \arr H$ is a splitting map, hence $\sigma h = id_G$. Notice that
the set
$$C = \{ \alpha < \omega_1 \;, H_\alpha \cap H^0 = H'_\alpha , G_\alpha \sigma
\subseteq H_\alpha\}$$

\noindent
is a cub. Since $S' \subseteq \omega_1$ is stationary, 
we find an $\alpha \in S' \cap C$
and also let $\sigma \upharpoonright G_\alpha = \sigma$, hence

\medskip

$(c)$ $\begin{cases}
\sigma : G_\alpha \arr H_\alpha \subset H^\alpha,\ P_\alpha (\sigma) =
\eta \; (\alpha) \mbox{ and }\\
\sigma \mbox{ is a splitting map of } h_\alpha : H^\alpha \arr G_\alpha.
\end{cases}$

\bigskip

\noindent
We also find some $\alpha < \beta \in C$.
The difficulty is that $G_{\alpha+1} \sigma \subseteq H_{\alpha+1}$ does not
follow, as in the case $S$ is not costationary. Hence we need the stronger Step-Lemma
(as usual).

\bigskip

\noindent
If $\eta (\alpha) = P (\sigma) = 0$, then $(H_{\eta \hat{~} 0}, h_{\eta \hat{~} 0})$
is part of the construction of 
$$0 \arr H^0 \arr H_{\alpha+1} \arr G_{\alpha+1} \arr 0$$
and $\sigma$ does not split over $(H_{\eta \hat {~} 0}, h_{\eta \hat {~} 0})$,
but $\sigma$ is a global splitting map, hence $\sigma$ splits at $\beta$ over
$(H_{\eta \hat{~} 0}, h_{\eta \hat {~} 0})$, a contradiction.\\
Necessarily $\eta (\alpha) = P (\sigma) = 1$ and by (a) $\sigma$ does not split
over $(H_{\eta \hat {~} 1}, h_{\eta \hat {~} 1})$, but this time
$(H_{\eta \hat {~} 1}, h_{\eta \hat {~} 1})$ was used in the construction of
$H \overset{h} \arr G$ and a contradiction follows. This shows that $\sigma$
is no splitting map, and $G$ is not a splitter.
$\hfill \square$\\

\noindent
R\"udiger G\"obel \\
Fachbereich 6, Mathematik und Informatik \\
Universit\"at Essen, 45117 Essen, Germany \\
{\small e--mail: R.Goebel@Uni-Essen.De}\\
and \\ 
Saharon Shelah \\
Department of Mathematics\\ 
Hebrew University, Jerusalem, Israel \\
and Rutgers University, Newbrunswick, NJ, U.S.A \\
{\small e-mail: Shelah@math.huji.ae.il}

\end{document}